\documentclass[a4paper,11pt]{article}
\usepackage{amsmath,amssymb,amsthm,enumerate}
\usepackage{xcolor}
\usepackage{graphicx}
\usepackage{epstopdf}
\usepackage{dsfont}
\usepackage{pdfpages}
\usepackage[left=2.6cm, right=2.6cm, top=2.5cm, bottom=2.5cm]{geometry}
\usepackage{url}
\usepackage[T1]{fontenc}
\usepackage[utf8]{inputenc}
\usepackage[ngerman,british,UKenglish,USenglish,american]{babel}
\usepackage{makeidx}
\usepackage{mathrsfs} 
\usepackage{setspace}
\usepackage{bbm}
\usepackage{bm}
\usepackage{relsize}
\usepackage[export]{adjustbox}
\usepackage{amssymb}
\usepackage{color}
\usepackage[toc,page]{appendix}

\usepackage[normalem]{ulem}
\usepackage{pdfpages}
\usepackage{caption}
\usepackage[ 
position=top]{subcaption}

\usepackage{hyperref}

\newcommand*{\cF}{\mathcal{F}}

\newcommand*{\cU}{\mathcal{U}}
\newcommand*{\cV}{\mathcal{V}}

\newcommand*{\R}{\mathds{R}}
\newcommand*{\E}{\mathds{E}}
\renewcommand*{\P}{\mathds{P}}

\newcommand{\be}{\begin{eqnarray*}}
	\newcommand{\ee}{\end{eqnarray*}}
\newcommand{\ben}{\begin{eqnarray}}
\newcommand{\een}{\end{eqnarray}}
\newcommand{\bi}{\begin{itemize}}
	\newcommand{\ei}{\end{itemize}}

\theoremstyle{plain}
\newtheorem{theo}{Theorem}[section]

\theoremstyle{definition}

\newtheorem{remark}[theo]{Remark}

\numberwithin{equation}{section}

\title{Ruin probability in a two-dimensional model with correlated Brownian motions}

\author{Peter Grandits\thanks{Peter Grandits,  TU Wien,
		Institute of Statistics and Mathematical Methods in Economics,
		Wiedner Hauptstr.\ 8 / E105-1\&-5 FAM,
		1040 Wien, Austria. \emph{Email:} pgrand@fam.tuwien.ac.at, \emph{Phone:} +43 1 58801 10512. } 
	\and Maike Klein\thanks{Maike Klein, TU Wien,
		Institute of Statistics and Mathematical Methods in Economics,
		Wiedner Haupt-\mbox{str.\ 8}\,/\,E105-1\&-5 FAM,
		1040 Wien, Austria. \emph{Email:} maike.klein@fam.tuwien.ac.at,\,\emph{Phone:} +43 1 58801 105170.\vspace*{0,2cm}\newline
				 Both authors gratefully acknowledge the support by the Austrian Science Fund (Fonds  zur Förderung der wissenschaftlichen Forschung) under grant P30864-N35.}
	}

\begin{document}

	\maketitle
\begin{abstract}
	
	We consider two insurance companies with endowment processes given by Brownian motions with drift. The firms can collaborate by transfer payments in order to maximize the probability that none of them goes bankrupt. 
	We show that pushing maximally the company with less endowment is the optimal strategy for the collaboration if the Brownian motions are correlated and the transfer rate can exceed the drift rates. 
	Moreover, we obtain an explicit formula for the minimal ruin probability in case of perfectly positively correlated Brownian motions where we also allow for different diffusion coefficients.
	
	\end{abstract}

\begin{center}\footnotesize
	\begin{tabular}{l@{ : }p{8.5cm}}
		{\itshape 2020 MSC} & Primary 49J20, 35R35, Secondary 91B70.\\
		{\itshape Keywords} & Ruin probabilities, Optimal control problem,  Collaboration, Two-dimensional Brownian motion, Correlated Brownian motions.
	\end{tabular}
\end{center}
\section*{Introduction}

We focus on two insurance companies whose endowment processes are given by correlated Brownian motions with drift. The common aim of the firms is to minimize the probability that at least one of the endowment processes falls below zero and, thus,  they collaborate by transfer payments.
These payments are assumed to be absolutely continuous with respect to the Lebesgue measure and to be bounded, but they can exceed the drift rates so that a company can be faced with a negative drift rate. The problem can be interpreted as an optimal control problem which consists in minimizing the probability that the two-dimensional endowment process leaves the positive quadrant and in identifying the optimal transfer payments. 

If at some point in time one firm has high endowment and the endowment of the other firm is close to zero then it seems reasonable that the latter is maximally supported and obtains the whole available drift rate. Using this so-called push-bottom strategy turns out to be optimal no matter how big the difference between the endowment processes is: The firm with less endowment receives the maximal drift.  To show this result we use a comparison principle for stochastic differential equations (SDEs) from \cite{IkedaWatanabe}. 

If the Brownian motions are perfectly positively correlated we derive a closed formula for the value function, because only one Brownian motion is involved and we can rewrite the ruin probability in terms for which explicit formulas are available. The arguments also apply to the case where the endowment processes have different diffusion coefficients and the Brownian motions are perfectly positively correlated. The value function turns out to be a classical solution of the corresponding Hamilton-Jacobi-Bellman equation.

This control problem was first studied by McKean and Shepp in \cite{McKeanShepp} for independent Brownian motions and the value function was derived in the case that the transfer payments are at most as high as the drift rates. 

In case that each company keeps a given minimal positive drift rate, the value function and the gain of collaboration for independent Brownian motions are obtained in \cite{Grandits2019_ruin_problem} by constructing a classical solution of the associated Hamilton-Jacobi-Bellman equation. 

Although the ruin probability is one of the most important evaluation criteria for insurance companies, there are only few articles dealing with two or more companies. For an overview of the one-dimensional case consult \cite{AsmussenAlbrecher}.
In \cite{ChanYangZhang2003} a two-dimensional model is analyzed and simple bounds for the ruin probability are obtained by using results from the one-dimensional case. 
The Laplace transform in the initial endowments of the probability that at least one of the two companies is ruined in finite time is derived in \cite{AvramPalmowskiPistorius2008ruin}.
Collamore \cite{Collamore1996} investigates the probability that a $d$-dimensional discrete process hits a $d$-dimensional set $A$ and obtains some large deviation results. For specific choices of the set $A$ the hitting probability can be seen as a ruin probability. 
The asymptotic behavior of the ruin probability if the initial endowments both tend to infinity under a light-tails assumption on the claim size distribution is analyzed in~\cite{AvramPalmowskiPistorius2008exit}. 

Let us emphasize that also for the maximal expected aggregated dividend payments, which is  another main evaluation criteria for insurance companies, the literature in the multidimensional setting  is scarce.  The optimal collaboration for maximizing the total dividend payments of two companies in different models is analyzed for example in  \cite{AlbrecherAzcueMuler, GerberShiu2006,  Grandits2019_dividend_problem, GuSteffensenZheng2018}.

\medskip
The paper is organized as follows. In Section \ref{sec:model} we introduce our model. We derive the optimal strategy for the transfer payments in order to minimize the ruin probability
in Section \ref{sec:optimal_strategy}. In Section \ref{sec:rho1_simple}  we focus on perfectly positively correlated Brownian motions and compute the value function explicitly. The same arguments are extended to a model with different diffusion coefficients in Section \ref{sec:rho1_ext}. Finally, we rewrite the minimal ruin probability for perfectly negatively correlated Brownian motions in terms of the hitting probability of a reflected Brownian motion with drift in Section \ref{sec:rho-1}.


\section{Model}\label{sec:model}
The endowment processes of the two companies are described by the stochastic processes 
\begin{align*}
X^x_t&=x+\mu_1 t +W_t+\int_0^t c_s\,  ds, \\
Y^y_t&=y+\mu_2 t+ \widehat{W}_t-\int_0^t c_s\, ds, 
\end{align*}
where $x,y>0$ denote the initial endowments, $\mu_1, \mu_2>0 $ are constant cash rates, e.g.\ premium rates, $\widehat{W}_t=\rho W_t+\sqrt{1-\rho^2}\,W^{(2)}_t$, $t\in[0,\infty)$, with $\big(W_t,W^{(2)}_t\big)_{t\in[0,\infty)}$ a two-dimensional Brownian motion. Here $\rho\in[-1,1]$ is the correlation coefficient of $(W_t)_{t\in[0,\infty)}$ and $(\widehat{W }_t)_{t\in[0,\infty)}$. The drift rates $(c_s)_{s\in[0,\infty)}$ can be interpreted as transfer payments from one company to the other one. More precisely, if $c_s>0$, then the first company obtains payments from the second company  at time $s$ and if $c_s<0$ it is vice versa. We say that the companies do not collaborate if $c_s=0$ for all $ s\in [0,\infty)$. 
We assume that the transfer payments are bounded in such a way that the total drift rate of each company is bounded below by $-\delta$, thus, 
\begin{align*}
c_s\in[-\mu_1-\delta,\,\mu_2+\delta]
\end{align*}
for some $\delta>-\min\{\mu_1,\mu_2\}$.

Introducing the control process $u_s:= \mu_1+c_s$ the endowment processes are given by
\begin{align}
\begin{split}\label{model}
X^{x,u}_t&=x+W_t+\int_0^t u_s\,  ds, \\
Y^{y,u}_t&=y+ \widehat{W}_t+\int_0^t\left(\bar\mu- u_s\right) ds, 
\end{split}
\end{align}
where $u_s\in[-\delta,\, \bar\mu+\delta]$ and $\bar\mu=\mu_1+\mu_2$. 

We aim at maximizing the probability that both firms survive forever. For this purpose denote by
\begin{align*}
\tau_X(x;u)&=\inf\{t\in[0,\infty)\colon X^{x,u}_t\leq 0\},\\
\tau_Y(y;u)&=\inf\{t\in[0,\infty)\colon Y^{y,u}_t\leq 0\}
\end{align*}
the ruin times of the first and second company, respectively, when the control $u$ is used. Let
\begin{align*}
\tau(x,y;u)=\tau_X(x;u)\wedge \tau_Y(y;u).
\end{align*}
Our target functional is then given by 
\begin{align*}
J(x,y;u)=\P[\tau(x,y;u)=\infty]
\end{align*}
and the value function is
\begin{align}\label{basic_pb}
V(x,y)=\sup_{u\in \cU} J(x,y;u), 
\end{align}
where $\,\cU$ denotes the set of all admissible controls. More precisely, $\cU$ is the set of all progressively measurable processes $(u_s)_{s\in[0,\infty)}$ with respect to the filtration generated by $\big(W_t, \widehat{W}_t\big)_{t\in[0,\infty)}$ satisfying $u_s\in[-\delta,\, \bar\mu+\delta]$, $s\in[0,\infty)$. 
 
 For $\bar\mu=1$ and $\rho=0$ we obtain the same model as in the paper by McKean and Shepp \cite{McKeanShepp}; for $\rho=0$ and $\delta\in(-\min\{\mu_1,\mu_2\},0)$ we are in the setting of \cite{Grandits2019_ruin_problem}. 
	
	The Hamilton-Jacobi-Bellman equation and the boundary conditions for the optimal control problem \eqref{basic_pb} are given by
	\begin{align*}
\left(	\frac 12 V_{xx}+\frac12 V_{yy}+\rho V_{xy}+(\bar\mu+\delta)\max\{ V_x,V_y\}-\delta\min\{V_x,V_y\}\right)(x,y)&=0\  \ \text{on $(0,\infty)\times(0,\infty)$,} \\[0.15cm]
	V(x,0)=V(0,y)&=0, \\[0.3cm]
	\lim_{x\to\infty} V(x,y)&= 1-\exp\big(-2(\bar\mu+\delta)y\big), \\[0.2cm]
		\lim_{y\to\infty} V(x,y)&= 1-\exp\big(-2(\bar\mu+\delta)x\big).
		\end{align*}

\begin{remark}
We can interpret \eqref{basic_pb}  as a tax policy problem, where the state can influence the endowment of the companies by imposing some kind of taxes, see \cite{McKeanShepp}. 
\end{remark}	
	
\section{The Optimal Strategy for the Transfer Payments}\label{sec:optimal_strategy}
For deriving the optimal strategy in the control problem \eqref{basic_pb} we use a comparison theorem for solutions of stochastic differential equations. We focus on the case $|\rho|<1$ because it is more involved and the arguments simplify for $|\rho|=1$ and, thus, are omitted. 

First, consider the transformation
\begin{align*}
Z^{(1)}_t= X_t^{x,u}+Y_t^{y,u}&= x+y+\bar\mu t+(1+\rho)W_t+\sqrt{1-\rho^2}\,W_t^{(2)}, \\[0.3cm]
Z^{(2), v}_t= Y_t^{y,u}-X_t^{x,u}&= y-x+\int_0^t (\bar\mu-2u_s)\,ds +(\rho-1)W_t+\sqrt{1-\rho^2}\,W_t^{(2)}\\[0.2cm]
&= y-x+\int_0^t v_s\, ds+(\rho-1)W_t+\sqrt{1-\rho^2}\,W_t^{(2)}, 
\end{align*}
where $v_s:=\bar\mu-2u_s$ and $v_s\in[-(\bar\mu+2\delta),\bar\mu+2\delta]$, $s\in[0,\infty)$. Observe that $Z^{(1)}$ does not depend on the control $v=(v_s)_{s\in[0,\infty)}$. 
Furthermore, it holds that 
\begin{align*}
\big\langle (1+\rho) W+\sqrt{1-\rho^2}\, W^{(2)}\big\rangle_t&= 2(1+\rho)t, \\
\big\langle (\rho-1) W+\sqrt{1-\rho^2}\, W^{(2)}\big\rangle_t&= 2(1-\rho)t, 
\end{align*}
and 
\begin{align*}
\left\langle (1+\rho) W+\sqrt{1-\rho^2}\, W^{(2)},  (\rho-1) W+\sqrt{1-\rho^2}\, W^{(2)}\right\rangle_{\!t}= \big((1+\rho)(\rho-1)+1-\rho^2\big)t=0. 
\end{align*}
We rescale the diffusion parts of the processes $Z^{(1)}$ and $Z^{(2),v}$ and define
\begin{align*}
B^{(1)}_t&=\sqrt{\frac{1+\rho}{2}}\, W_t+\sqrt{\frac{1-\rho}{2}}\, W^{(2)}_t,\qquad
B^{(2)}_t=-\sqrt{\frac{1-\rho}{2}}\, W_t+\sqrt{\frac{1+\rho}{2}}\, W^{(2)}_t.
\end{align*}
Hence, we conclude that $\big(B^{(1)}_t)_{t\in [0,\infty)}$, $\big(B^{(2)}_t\big)_{t\in [0,\infty)}$ are independent Brownian motions and that 

\begin{align*}
Z^{(1)}_t &= x+y+\bar\mu t+\sqrt{2(1+\rho)}\,B^{(1)}_t, \\[0.2cm]
Z^{(2), v}_t&= y-x+\int_0^t v_s\, ds+\sqrt{2(1-\rho)}\,B_t^{(2)},
\end{align*}
where  $(v_s)_{s\in[0,\infty)}\in\cV$. Here $\cV$ denotes the set of all progressively measurable processes $(\tilde v_s)_{s\in[0,\infty)}$
with respect to the filtration $(\cF_t)_{t\in[0,\infty)}$, which is generated by $\big(\sqrt{2(1+\rho)}B^{(1)}_t\big)_{t\in[0,\infty)}$ and  $\big(\sqrt{2(1-\rho)}B^{(2)}_t\big)_{t\in[0,\infty)}$, and $(\tilde v_s)_{s\in[0,\infty)}$ satisfies $\tilde v_s~\in~[-(\bar\mu+2\delta),\bar\mu+2\delta]$, $s\in[0,\infty)$. 
 Note that $(\cF_t)_{t\in[0,\infty)}$ and the  filtration generated by $\big(W_t, \widehat{W}_t\big)_{t\in[0,\infty)}$ coincide.

 So far all arguments hold for $\rho\in[-1,1]$.  Now let $|\rho|<1$. Then $(\cF_t)_{\in[0,\infty)}$ also coincides with the filtration generated by the Brownian motion $\big(B_t^{(1)}, B_t^{(2)}\big)_{t\in[0,\infty)}$. In particular, there exists a measurable function 
   \begin{align*}
v\colon [0,\infty)\times C\big([0,\infty);\R\big)\times C\big([0,\infty);\R\big)\to [-(\bar\mu+2\delta), \bar\mu+2\delta]
  \end{align*}
  such that for every $s\in[0,\infty)$
\begin{align*}
v_s= v\left(s,\big(B^{(1)}_{r\wedge s}\big)_{r\in [0,\infty)},\big(B^{(2)}_{r\wedge s}\big)_{r\in [0,\infty)}\right). 
\end{align*}
We rewrite $\tau(x,y;u)$ in terms of $Z^{(1)}$ and $Z^{(2),\, \bar\mu-2u}$. More precisely, 
\begin{align*}
\tau(x,y;u)&
=\inf\left\{t\in[0,\infty)\colon \Big|Z_t^{(2),\,  \bar \mu-2u}\Big|\geq Z_t^{(1)}\right\}=:\tau^{\bar\mu-2u}.
\end{align*}
Hence,
\begin{align}\label{transformed_pb}
V(x,y)=\sup_{v\in \cV}\P[\tau^v=\infty]
\end{align}
and for every optimal $v^*\in\cV$ in \eqref{transformed_pb} we obtain an optimal $u^*\in \cU$ for \eqref{basic_pb} by setting $u^*=(\bar\mu-v^*)/2$.

To characterize an optimal control $v^*$ for \eqref{transformed_pb} denote by $\cF^{(1)}$ the filtration generated by $\big(B^{(1)}_t\big)_{t\in[0,\infty)}$, i.e.
{\setlength{\abovedisplayskip}{0pt}	
	\setlength{\belowdisplayskip}{3pt}
	\begin{align*}
\cF^{(1)}=\sigma\left(B_s^{(1)}\colon s\in[0,\infty)\right),
\end{align*}}
\!\!and observe that 
\begin{align*}
\P[\tau^{v}=\infty]&=\E\left[\E\left[\mathds{1}_{\left\{\big|Z^{(2), v}_t\big|\,<\,Z^{(1)}_t\ \forall\, t\,\in\,[0,\infty)\right\}}\Big|\cF^{(1)}\right] \right]\\
&=\E\left[\E\left[\mathds{1}_{\left\{\big|y-x+\mathlarger{\int}_0^t v\big(s,\big(B_{r\wedge s}^{(1)}\big)_{r\in[0,\infty)}, \big(B_{r\wedge s}^{(2)}\big)_{r\in[0,\infty)}\big)ds\,+\sqrt{2(1-\rho)}\,B_t^{(2)}\big|\,<\,Z^{(1)}_t\ \forall \,t\,\in\,[0,\infty)\right\}}\Big|\cF^{(1)}
\right]\right]. 
\end{align*}
Note that $\big(B_t^{(1)}\big)_{t\in[0,\infty)}$ and $\big(Z^{(1)}_t \big)_{t\in[0,\infty)}$ are measurable with respect to $\cF^{(1)}$ and that $\big(B^{(2)}_t\big)_{t\in[0,\infty)}$ is independent of $\cF^{(1)}$.
Now choose a realization of $B^{(1)}$. In particular, $B_t^{(1)}:= f(t)$ and $Z^{(1)}_t:= g(t)$ are fixed. Furthermore, let 
\begin{align*}
\hat{v}\big(s,\big(B^{(2)}_{r\wedge s})_{r\in[0,\infty)}\big)\big)= v\big(s,(f(r\wedge s))_{r\in[0,\infty)}, (B_{r\wedge s}^{(2)})_{r\in[0,\infty)}\big). 
\end{align*}
If the optimal strategy $\hat{v}^*$ for maximizing 
\begin{align}\label{aux_max}
\E&\left[\mathds{1}_{\left\{\big|y-x+\mathlarger{\int}_0^t \hat{v}\big(s, (B_{r\wedge s}^{(2)})_{r\in[0,\infty)}\big)ds\,+\sqrt{2(1-\rho)}\,B_t^{(2)}\big|\,<\,g(t)\ \forall \,t\,\in\,[0,\infty)\right\}}\right]
\end{align}
is independent of $f$ and $g$, it is also optimal for maximizing $\P[\tau^v=\infty]
$ over all $v\in\cV$.  
To identify $\hat{v}^*$ one can extend the arguments from Theorem 2.1 in \cite{IkedaWatanabe}  and its corollary to the case, where the control does not take values in $[-1,1]$ but in $[-(\bar\mu+2\delta), \bar\mu+2\delta]$ and the diffusion coefficient equals $\sqrt{2(1-\rho)}$ instead of 1. Observe that the required filtration is given by $(\cF_t)_{t\in[0,\infty)}$ and  that $\big(B^{(2)}_t\big)_{t\in[0,\infty)}$ is an $(\cF_t)_{t\in[0,\infty)}$-Brownian motion.  Then for all $\hat{v}$ we have 
\begin{align*}
\P\left[\left|y-x+\int_0^t \hat{v}\big(s, (B_{r\wedge s}^{(2)})_{r\in[0,\infty)}\big)ds+\sqrt{2(1-\rho)}\,B_t^{(2)}\right|<g(t)\ \forall\, t\in[0,\infty)\right]&\\[0.2cm]
&\hspace*{-2cm}\leq \P\left[Z^*_t<g(t)\ \forall\, t\in[0,\infty)\right],
\end{align*}
where $(Z^*_t)_{t\in[0,\infty)}$ satisfies
\begin{align}\label{SDE}
dZ_t^*= -(\bar\mu+2\delta)\, \text{sign}(Z^*_t)\,dt+ \sqrt{2(1-\rho)}\,dB_t^{(2)}, \quad Z_0^*=y-x. 
\end{align}
with
\begin{align*}
	\text{sign}(z)=\begin{cases}
		1, \quad & z>0, \\
		0, &z=0,\\
		-1, & z<0. 
	\end{cases}
\end{align*}
The existence of a strong solution of \eqref{SDE} which is strongly unique follows from Theorem~1 in~\cite{veretennikov1983}.
In particular, there exists a measurable function $F$ such that 
\begin{align*}
Z^*_t=F\big(\big(B_{s\wedge t}^{(2)}\big)_{s\in[0,\infty)}\big), \quad t\in[0,\infty),
\end{align*}
 and, thus, the drift of $Z^*_t$ is a  measurable function of the Brownian motion $B^{(2)}$ up to time $t\in[0,\infty)$. Hence, the optimal control $\hat{v}^*$ for maximizing \eqref{aux_max} is given by
\begin{align*}
\hat{v}^*\big(s,(B_{r\wedge s}^{(2)})_{r\in[0,\infty)}\big)=-(\bar\mu+2\delta)\, \text{sign}\!\left(F\left((B_{r\wedge s}^{(2)})_{r\in[0,\infty)}\right)\right), \quad s\in[0,\infty), 
\end{align*}
and, in particular, it is independent of $f$ and $g$. Therefore, 
\begin{align*}
v_s^*=-(\bar\mu+2\delta)\,\text{sign}\big(Z_s^{(2), v^*}\big), \quad s\in[0,\infty),
\end{align*}
is optimal for  \eqref{transformed_pb}. By setting $u^*_s=\frac{\bar\mu-v^*}{2}$ we obtain an optimal control for \eqref{basic_pb}. 

For $|\rho|=1$ one also uses the arguments from \cite{IkedaWatanabe} and obtains the same optimal strategy.

\begin{remark}\label{change_on_nullset}
	Observe that one can change the definition of the optimal strategy $u^*$  on the set $\{X^{x,u^*}_t=Y^{y,u^*}_t\}=\{Z^*_t=0\}$, $t\in[0,\infty)$, and then obtains an indistinguishable process, because with probability one the set $\{t\in[0,\infty)\colon Z^*_t=0\}$, where $(Z^*_t)_{t\in[0,\infty)}$ is a solution of~\eqref{SDE}, has Lebesgue measure zero, for details see Appendix C in \cite{Benes1976}.
\end{remark}

We now summarize the result in the following theorem. 

\begin{theo}\label{optimal}
	Let $\rho\in[-1,1]$. The optimal drift rate $(u_s^*)_{s\in[0,\infty)}$ in \eqref{basic_pb} for minimizing the ruin probability is given by
	\begin{align*}
	u_s^*=u_s^*\left(X_s^{x,u^*}, Y_s^{y,u^*}\right)=(\bar\mu+\delta)\mathds{1}_{\left\{X_s^{x,u^*}\!\leq\, Y_s^{y,u^*}\right\}}
		-\delta\mathds{1}_{\left\{X_s^{x,u^*}\!>\,Y_s^{y,u^*}\right\}}, \quad\ s\in[0,\infty).
	\end{align*}
	Moreover, the optimal transfer rate is
	$c_s^*= u^*_s-\mu_1$.
\end{theo}

Theorem \ref{optimal} implies that the company with less endowment is as much supported as possible, i.e.\ this company receives a drift of $\bar\mu+\delta$ until the endowment processes of the two companies are equal: The push-bottom strategy is optimal.

\begin{remark}
In \cite{FernholzIchibaKaratzasProkaj} and \cite{FernholzichibaKaratzas} the authors analyze two diffusion processes where the drift and diffusion coefficients are rank-dependent and the Brownian motions are independent. The leader obtains a negative drift coefficient and the laggard a positive one.
For the isotropic case, i.e.\ for the same diffusion coefficients, this corresponds to the optimal controlled processes $X^{x,u^*}$ and $Y^{y,u^*}$ for $\rho=0$ in our setting.
\end{remark}

\begin{remark}
	For independent Brownian motions, i.e.\ $\rho=0$, the value function can be computed explicitly. 
	McKean and Shepp \cite{McKeanShepp} show that for $\delta=0$ and $\bar\mu=1$ the value function is given by
	\begin{align*}
	V(x,y)=1-e^{-2\min\{x,y\}}-2\min\{x,y\}\,e^{-x-y}. 
	\end{align*}
	For  $\delta\in(-\min\{\mu_1,\mu_2\},0)$ Grandits obtains in Theorem 4.1 of  \cite{Grandits2019_ruin_problem}  the following value function
	\begin{align*}
	V(x,y)=1-e^{-(\bar\mu+\delta)\min\{x,y\}}-\frac{\bar\mu+\delta}{\delta}e^{-\bar\mu(x+y)}\left(1-e^{-2\delta\min\{x,y\}}\right).
	\end{align*}
	\end{remark}

\section{Perfectly Positive Correlation: $\rho=1$}\label{sec:rho1}
In this section we consider perfectly positively correlated Brownian motions, i.e.\ $\rho=1$, which implies $\widehat{W}_t=W_t$, and derive an explicit formula for the value function \eqref{basic_pb}. In Section \ref{sec:rho1_simple} we deal with the simplest case where $\delta=0$ and $\bar\mu=1$. 
For $\rho=1$ the same arguments can be used for endowment processes $X$ and $Y$ having different diffusion coefficients. Hence,
 we extend our model and state the value function for different diffusion coefficients $\sigma_1>0$ and $\sigma_2>0$ in Section \ref{subsec:diff_coeff}. In both sections we also compute the gain of collaboration.

\subsection{Deriving the Value Function for $\delta=0$ and $\bar\mu=1$}\label{sec:rho1_simple} 
We now derive the value function for $\delta=0$ and $\bar\mu=1$. The same arguments extend to the case $\delta >-\min\{\mu_1,\mu_2\}$ with $\mu_1,\mu_2>0$ but lead to more complicated terms. Therefore, we first focus on this simple case. 

Let $\delta=0$, $\bar\mu=1$.
As in Section \ref{sec:optimal_strategy} consider the processes
\begin{align*}
Z^{(1)}_t&= X_t^{x,u}+Y_t^{y,u}= x+y+2 W_t+ t , \\
Z^{(2), v}_t&= Y_t^{y,u}-X_t^{x,u}= y-x+\int_0^t (1-2u_s)\,ds= y-x+\int_0^t v_s\, ds, 
\end{align*}
where $v_s:=1-2u_s$ and $v_s\in [-1,1]$, $s\in[0,\infty)$. Recall that $Z^{(1)}$ does not depend on the control $v=(v_s)_{s\in[0,\infty)}$,
 \begin{align*}
 \tau(x,y;u)=\inf\left\{t\in[0,\infty)\colon  Z_t^{(1)}\leq \left| Z_t^{(2), 1-2u}\right| \right\}=:\tau^{1-2u}. 
 \end{align*}
and that 
\begin{align*}
V(x,y)=\sup_{v\in \cV}\P[\tau^v=\infty].
\end{align*}
Here $\cV$ is the set of all progressively measurable processes $(v_s)_{s\in[0,\infty)}$ with respect to the filtration generated by $(W_t)_{t\in[0,\infty)}$ satisfying $v_s~\in~[-1,\,1]$ for all $s\in[0,\infty)$.

Since $Z^{(1)}$ is a Brownian motion with drift and it is not controllable by $v$ and $Z^{(2),v}$ does not depend on the Brownian motion $(W_t)_{t\in[0,\infty)}$, the best strategy is to control $Z^{(2),v}$ in such a way that it is going to zero with the highest possible rate and then let $Z^{(2),v}$  stay zero. In particular, 
\begin{align*}
v_s^*=-\,\text{sign}\! \left(Z_s^{(2),v^*}\right), \qquad s\in[0,\infty),
\end{align*}
is optimal.
Note that the corresponding optimal strategy $u^*$ in \eqref{basic_pb} is given by $u_s^*=\mathds{1}_{\{X_s^{x,u^*}\!\leq\, Y^{y,u^*}_s\}}$, which we have already seen in Theorem \ref{optimal}; also recall Remark \ref{change_on_nullset}.

\begin{figure}
	\begin{center}
		\includegraphics[width=0.7\textwidth]{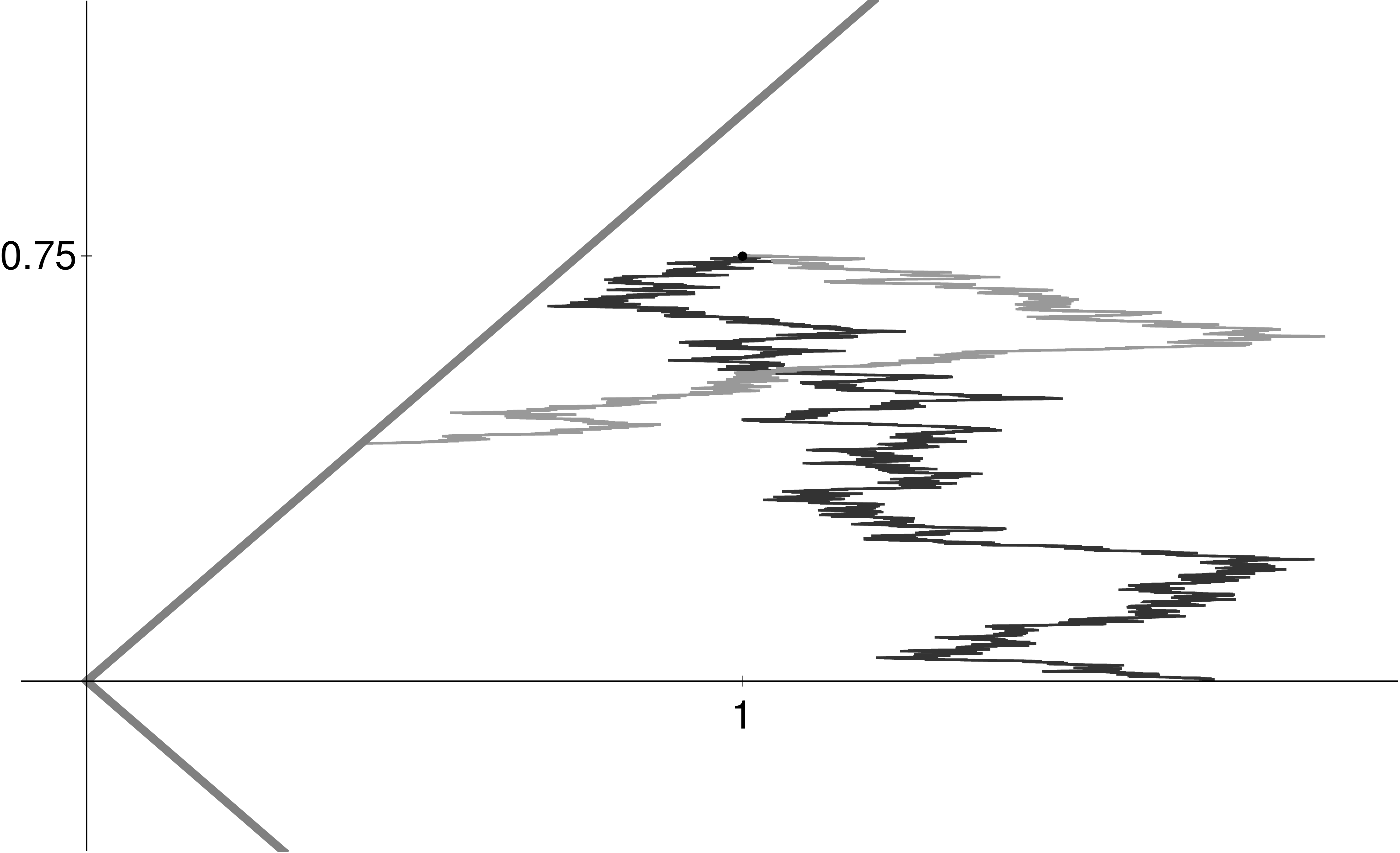}
	\end{center}
	\caption{Two trajectories of $(Z^{(1)}, Z^{(2),v^*})$.  In case of the light gray trajectory ruin occurs before $Z^{(2), v^*}$ becomes zero, i.e.\ before time $y-x$. In the other case ruin has not occurred until time $y-x$.  }\label{fig:trajectories_stop}
\end{figure}

Now assume that $y> x$. Then the optimal controlled process $Z^{(2),v^*}$ is given by
\begin{align*}
Z_t^{(2),v^ *}= \max\{y-x-t, 0\}
=\begin{cases} 
y-x-t, \quad  & t\leq y-x,\\
0, & t\geq y-x. 
\end{cases}
\end{align*}
For computing the value function, we focus on $\P[\tau^{v^*}<\infty]$. On the set $\{\tau^{v^*}<\infty \}$ the stopping time $\tau^{v^*}$ either occurs before the process $Z^{(2),v^*}$ becomes zero or afterwards. For two possible trajectories see Figure \ref{fig:trajectories_stop}.
 Hence, 
\begin{align*}
\P[\tau^{v^*}\!<\infty]&= \P[ \tau^{v^*}\!\leq y-x]+ \P[y-x<\tau^{v^*}\!<\infty]\\[0.3cm]
&=\P\left[  Z_t^{(1)}\leq Z_t^{(2),v^*}\!\text{ for some $t\in[0,y-x]$}\right]+\int_{(0,\infty)}\!\P\left[y-x<\tau^{v^*}\!\!<\infty,\, Z^{(1)}_{y-x}\in dw\right]\!.
\end{align*}
We compute the two summands separately. 
\begin{align}
\begin{split}\label{p1}
\P\left[  Z_t^{(1)}\leq Z_t^{(2),v^*}\! \text{ for some $t\in[0,y-x]$}\right]&=\P\left[  \inf_{0\leq t\leq y-x} \big\{W_t+t\big\}\leq -x\right]\\[0.3cm]
& = 1-\Phi\left(\frac{y}{\sqrt{y-x}}\right)+\exp\left(-2x\right)\Phi\left(\frac{y-2x}{\sqrt{y-x}}\right),
\end{split}
\end{align}
where $\Phi$ denotes the cumulative distribution function of a standard normal distribution. The last equality follows from Formula 1.2.4 in Part II, Chapter 2 of \cite{Borodin}. 

For the second summand we obtain
\begin{align}
\int_{(0,\infty)}&\P\left[y-x<\tau^{v^*}\!<\infty,\, Z^{(1)}_{y-x}\in dw\right]\notag\\
&= \int_{(0,\infty)} \P\left[\tau^{v^*}\!<\infty\,\Big|\, \tau^{v^*}\!>y-x,\, Z^{(1)}_{y-x}=w\right]\, \P\left[\tau^{v^*}\!>y-x,\, Z^{(1)}_{y-x}\in dw\right]\notag\\[0.3cm]
&= \int_{(0,\infty)} \P\left[Z^{(1)}_t \leq 0 \,\text{ for some $t>y-x$}\,\Big|\, \tau^{v^*}\!>y-x, Z^{(1)}_{y-x}=w\right]\, \P\left[\tau^{v^*}\!>y-x,\, Z^{(1)}_{y-x}\in dw\right]\notag\\[0.3cm]
&= \int_{(0,\infty)} \P\left[\,\inf_{t\geq 0} \Big\{w+t+2\big(W_{y-x+t}-W_{y-x}\big)\Big\}\leq 0\right]\, \P\left[\tau^{v^*}\!>y-x,\, Z^{(1)}_{y-x}\in dw\right]\notag\\[0.3cm]
&=\int_{(0,\infty)} \exp\left(-\frac w2\right)\, \P\left[\tau^{v^*}\!>y-x,\, Z^{(1)}_{y-x}\in dw\right],\label{p21}
\end{align}
where the last equality follows from Formula 1.2.4 (1) in Part II, Chapter 2 of \cite{Borodin}. For the remaining probability we have
\begin{align}
\P[\tau^{v^*}\!>y-x&,\, Z^{(1)}_{y-x}\in dw]= \P\left[Z^{(1)}_t>Z^{(2),v^*}_t \!\text{ for all $t\in[0,y-x]$},\, Z^{(1)}_{y-x}\in dw\right]\notag\\[0.3cm]
&= \P\left[\inf_{0\leq t\leq y-x} \{W_t+t\}>-x,\, y+W_{y-x}\in \frac{dw}{2}\right]\notag\\[0.3cm]
&= \P\left[y+W_{y-x}\in \frac{dw}{2}\right]- \P\left[\inf_{0\leq t\leq y-x} \{W_t+t\}\leq -x,\, W_{y-x}+y-x\in \frac{dw}{2}-x\right]\notag\\[0.3cm]
&=\frac{1}{2\sqrt{2\pi(y-x)}}\left[\exp\left(-\frac{(w-2y)^2}{8(y-x)}\right) -\exp\left(  \frac{w-x-y}{2}-\frac{(w+2x)^2}{8(y-x)}\right)\right]dw,\label{p22} 
\end{align}
where we used Formula 1.2.8 in Part II, Chapter 2 of \cite{Borodin} in the last equality.
Hence, combining  \eqref{p21} and \eqref{p22} yields
\begin{align}\label{p2}
\int_{(0,\infty)}&\P\left[y-x<\tau^{v^*}\!<\infty,\, Z^{(1)}_{y-x}\in dw\right]= \exp\left(-\frac{x+y}{2}\right)\left[2\Phi\left(\frac{x}{\sqrt{y-x}}\right)-1\right]. 
\end{align}
Therefore, we conclude from \eqref{p1} and \eqref{p2} that for $y>x$ it holds that 
\begin{align*}
V(x,y)=\Phi\left(\frac{y}{\sqrt{y-x}}\right)-\exp\left(-2x\right)\Phi\left(\frac{y-2x}{\sqrt{y-x}}\right)-\exp\left(-\frac{x+y}{2}\right)\left[2\Phi\left(\frac{x}{\sqrt{y-x}}\right)-1\right]. 
\end{align*}
For $y=x$ we have that $Z^{(2),v^*}_t\!\!=0$ for all $t\in [0,\infty)$. Therefore, 
\begin{align*}
V(x,x)= \P\left[Z^{(1)}_t>0 \text{ for all $t\in [0,\infty)$}\right]=1-\exp\left(-x\right). 
\end{align*}
For $y<x$ we use the symmetry of the problem and conclude that $V(x,y)=V(y,x)$. 

To summarize, we have shown the following result. 

\begin{theo}\label{value_function_simple}
For $\rho=1$, $\delta=0$ and $\bar\mu=1$ the value function \eqref{basic_pb} is given by
\begin{align*}
V(x,y)=\begin{cases}
\begin{aligned}
&1-\exp(-x)= 1-\exp(-y), && y=x,\\[0.5cm]
&\Phi\left(\frac{\max\{x,y\}}{\sqrt{|y-x|}}\right)
-\exp\left(-2\min\{x,y\}\right)\Phi\left(\frac{|y-x|-\min\{x,y\}}{\sqrt{|y-x|}}\right) \quad &&y\neq x. \\[0.2cm]
&\quad -\exp\left(-\frac{x+y}{2}\right)\left(2\Phi\left(\frac{\min\{x,y\}}{\sqrt{|y-x|}}\right)-1\right),&&
\end{aligned}
\end{cases}
\end{align*}
\end{theo}

\begin{figure}
	\begin{center}
		\includegraphics[width=0.55\textwidth]{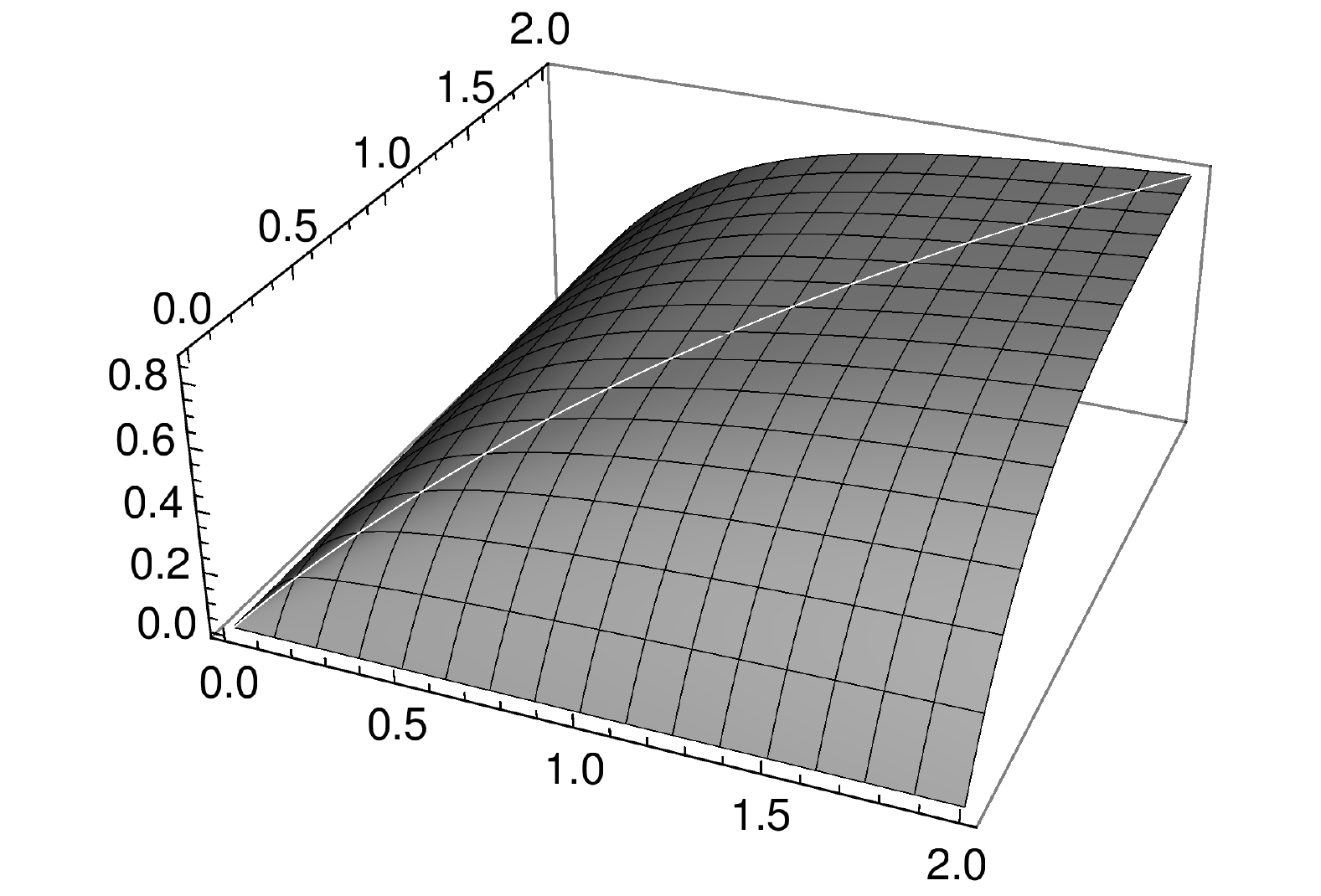}
		\caption{The value function $V$ for $\rho=1$, $\delta=0$ and $\bar\mu=1$.}\label{fig:value} 
	\end{center}
\end{figure}

Figure \ref{fig:value} depicts the value function of Theorem \ref{value_function_simple}. 

\begin{remark}
One can show that the function $V$ stated in Theorem \ref{value_function_simple} satisfies 
\begin{align*}
	V\in C^2\big((0, \infty)^2\big)\cap C\big([0,\infty)^2\big),
	\end{align*}
	 $V$ solves
\begin{align*}
\frac 12 V_{xx}(x,y)+\frac 12 V_{yy}(x,y)+V_{xy}(x,y)+\max\{V_x(x,y), V_y(x,y)\}&=0,\\[0.2cm]
V(x,0)=V(0,y)&=0, \\[0.25cm]
\lim_{x\to\infty} V(x,y)&= 1-\exp(-2y),\\[0.2cm]
\lim_{y\to\infty} V(x,y)&= 1-\exp(-2x),\\
\end{align*}
and
\begin{align*}
V_x(x,y)-V_y(x,y)=\begin{cases}
\begin{aligned}
&2\exp(-2x)\Phi\left(\frac{y-2x}{\sqrt{y-x}}\right)>0,\quad  && x <y,\\[0.4cm]
&-2\exp(-2y)\Phi\left(\frac{x-2y}{\sqrt{x-y}}\right)<0,\quad  && x >y,\\[0.4cm]
&0, &&x=y. 
\end{aligned}
\end{cases}
\end{align*}
Therefore, 
\begin{align*}
\left\{(x,y)\in(0,\infty)^2\colon V_x(x,y)>V_y(x,y)\right\}&=\left\{(x,y)\in(0,\infty)^2\colon x<y\right\},\\
\left\{(x,y)\in(0,\infty)^2\colon V_x(x,y)<V_y(x,y)\right\}&=\left\{(x,y)\in(0,\infty)^2\colon x>y\right\}, \\
\left\{(x,y)\in(0,\infty)^2\colon V_x(x,y)=V_y(x,y)\right\}&=\left\{(x,y)\in(0,\infty)^2\colon x=y\right\}.
\end{align*}

Hence, another possibility for proving Theorem~\ref{value_function_simple} is to use a classical verification theorem.  For applying this standard method the explicit formula of the value function has to be known in advance. The advantage of our approach is that it allows to  compute the optimal strategy and the value function directly. 
\end{remark}

We now compare the gain of collaboration. If the two firms do not collaborate, i.e.\ $u_s=\mu_1$ for all $s\in[0,\infty)$, then the survival probability of both firms is given by
\begin{align}
\P[X_t^{x,\mu_1}>0,\, Y_t^{y,\mu_1}>0 \text{ for all $t\in[0,\infty)$}]. \label{survival_prob}
\end{align}
For the cases $\mu_1=\mu_2$, $\mu_1>\mu_2$  with $x\geq y$ and $\mu_1<\mu_2$ with 
$ x\leq y$ the endowment of one company is for all $t\in [0,\infty)$ lower than the other company's endowment. Thus, \eqref{survival_prob} is just the survival probability of the firm with lower endowment and we have
\begin{align*}
\P[X_t^{x,\mu_1}>0,\, Y_t^{y,\mu_1}>0 \text{ for all $t\in[0,\infty)$}]=1-\exp\big(-2\min\left\{\mu_1 x, \mu_2 y\right\}\big).
\end{align*}
For ${\mu_1}>{\mu_2}$ with $x<y$ and ${\mu_1}<{\mu_2}$ with $ x>y$ the company with lower initial endowment has less endowment  until time $\frac{y-x}{\mu_1-\mu_2}$ and afterwards its endowment process is always larger than the process of the company with higher initial endowment. Therefore, it becomes more involved to compute \eqref{survival_prob}. We apply similar arguments as in the derivation of the value function in Theorem \ref{value_function_simple}, in particular, we use Formulas 1.2.4, 1.2.4 (1) and 1.2.8 in Part II, Chapter 2 of \cite{Borodin}.

{\captionsetup[figure]{skip=4pt}
	\begin{figure}
		\begin{subfigure}{0.48\textwidth}
			\includegraphics[width=\textwidth]{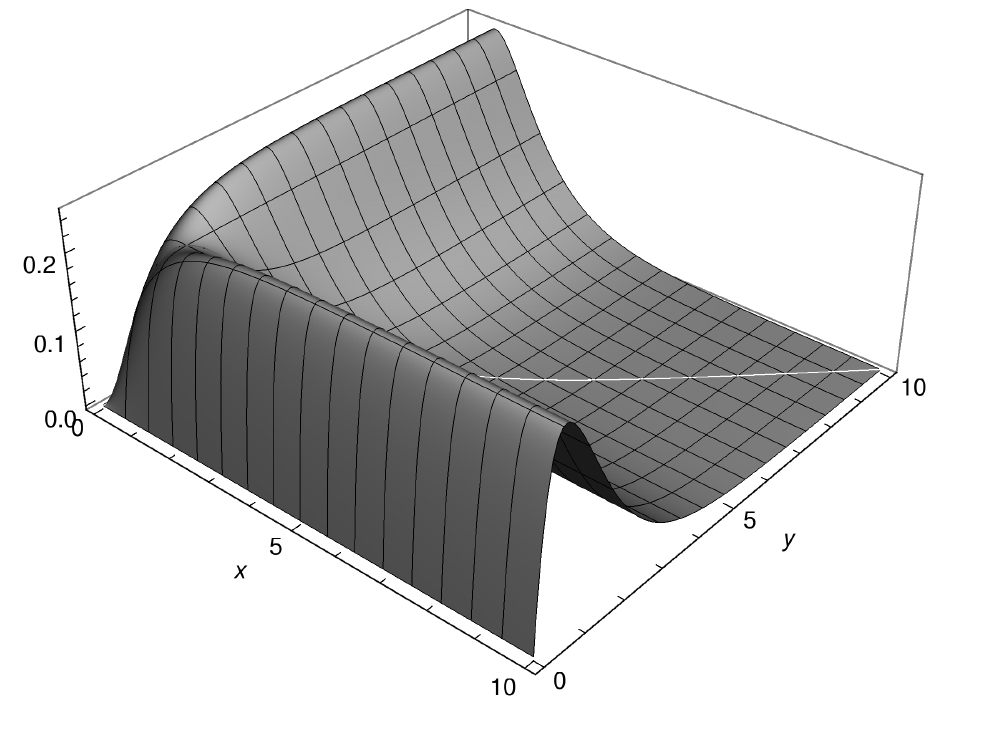}
			\subcaption{$\mu_1=\mu_2=\frac 12.$}
		\end{subfigure}
		\hspace{0.02\textwidth}
		\begin{subfigure}{0.48\textwidth}
			\includegraphics[width=\textwidth]{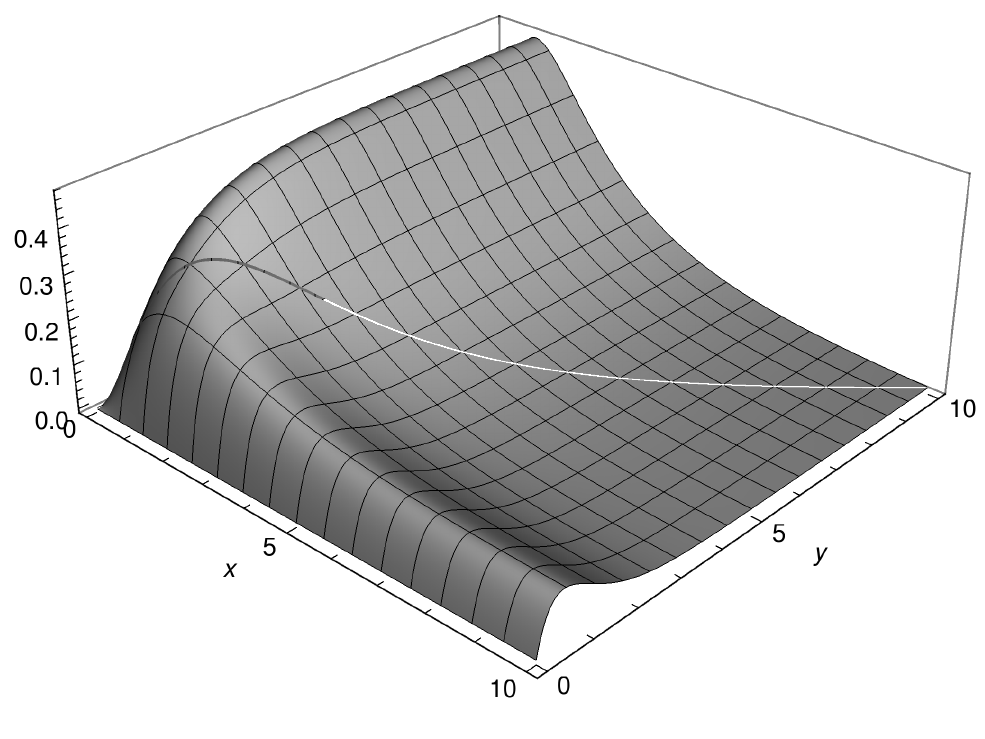}
			\subcaption{$\mu_1=\frac 14, \mu_2=\frac 34.$}
		\end{subfigure}
		\caption{The gain of collaboration for $\rho=1$, $\delta=0$, $\bar\mu=1$ and different drift rates $\mu_1$ and $\mu_2$. }\label{fig:gain}
\end{figure}}
Assume that $\mu_1>\mu_2 $ and $x<y$. Then it holds that
\begin{align*}
\P[&X_t^{x,\mu_1}>0,\, Y_t^{y,\mu_1}>0 \text{ for all $t\in[0,\infty)$}]\\[0.3cm]
&= 1-\P\left[\inf_{0\leq t\leq \frac{y-x}{\mu_1-\mu_2}}X_t^{x,\mu_1}\leq 0\right ]-\P\left[\inf_{0\leq t\leq \frac{y-x}{\mu_1-\mu_2}}X_t^{x,\mu_1}> 0, \inf_{ t\geq \frac{y-x}{\mu_1-\mu_2}}Y_t^{y,\mu_1}\leq 0\right]\\[0.3cm]
&=\Phi\left(\frac{\mu_1 y-\mu_2 x}{\sqrt{(y-x)(\mu_1-\mu_2)}}\right)
 +\exp\big((4\mu_2-2\mu_1) x-2\mu_2 y\big)\Phi\left(\frac{(3\mu_2-2\mu_1)x+(\mu_1-2\mu_2)y}{\sqrt{(y-x)(\mu_1-\mu_2)}}\right)
 \\[0.3cm]
 &\quad -\exp\left(-{2\mu_1 x}\right)\Phi\left(\frac{\mu_1 y+(\mu_2-2\mu_1)x}{\sqrt{(y-x)(\mu_1-\mu_2)}}\right)  
 -\exp\big(-{2\mu_2 y}\big)\Phi\left(\frac{\mu_2 x +(\mu_1-2\mu_2)y}{\sqrt{(y-x)(\mu_1-\mu_2)}}\right). 
\end{align*}
For ${\mu_1}<{\mu_2}$ with $ x>y$ change the role of $x$ and $y$ and the role of $\mu_1$ and $\mu_2$. 

The gain of collaboration is then given by $V(x,y)-\P[X_t^{x,\mu_1}>0,\, Y_t^{y,\mu_1}>0 \text{ for all $t\in[0,\infty)$}]$. See Figure \ref{fig:gain} for the gain of collaboration for different drift rates $\mu_1$ and $\mu_2$.  

\subsection{Different Diffusion Coefficients \label{sec:rho1_ext}
}\label{subsec:diff_coeff}
We now extend the model \eqref{model} and allow for different diffusion coefficients for $X$ and $Y$. More precisely, for $\sigma_1, \sigma_2>0$ let  
\begin{align}
\begin{split}\label{model_mod}
\widetilde{X}^{x,u}_t&=x+\sigma_1 W_t+\int_0^t u_s\,  ds, \\
\widetilde{Y}^{y,u}_t&=y+ \sigma_2 W_t+\int_0^t\left(\bar\mu- u_s\right) ds, 
\end{split}
\end{align}
where $u_s\in[-\delta\sigma_1,\, \bar\mu+\sigma_2\delta]$, $\delta>-\frac{\bar\mu}{\sigma_1+\sigma_2}$, $\bar\mu>0$. 
Thus, the relative drift rates $\frac{u_s}{\sigma_1}$ and $\frac{\bar\mu-u_s}{\sigma_2}$ are bounded below by $-\delta$.
 
Define $\widetilde{\tau}(x,y;\sigma_1, \sigma_2;u)=\tau_{\widetilde{X}}(x;u)\wedge\tau_{\widetilde{Y}}(y;u)$ and denote by 
\begin{align*}
V(x,y; \sigma_1, \sigma_2)=\sup_{u\in\widetilde\cU}\, \P[\,\widetilde{\tau}(x,y;\sigma_1, \sigma_2;u)=\infty]
\end{align*}
the value function in the extended model \eqref{model_mod}. Here $\widetilde\cU$ denotes the set of progressively measurable processes $(\tilde u_s)_{s\in[0,\infty)}$ with respect to the filtration generated by $(W_t)_{t\in[0,\infty)}$ such that $\tilde u_s\in [-\delta \sigma_1, \bar\mu + \delta\sigma_2]$ for all $s\in[0,\infty)$. Consider the transformation
\begin{align*}
\widetilde{Z}^{(1)}_t&= \widetilde{X}_t^{x,u}+\widetilde{Y}_t^{y,u}= x+y+(\sigma_1+\sigma_2)\,W_t+ \bar\mu t , \\[0.3cm]
\widetilde{Z}^{(2), v}_t&= \sigma_1 \widetilde{Y}_t^{y,u}-\sigma_2 \widetilde{X}_t^{x,u}= \sigma_1y-\sigma_2x+\int_0^t \big(\sigma_1\bar\mu-(\sigma_1+\sigma_2)u_s\big)\,ds= \sigma_1y-\sigma_2x+\int_0^t v_s\, ds, 
\end{align*}
where $v_s:=\sigma_1\bar\mu-(\sigma_1+\sigma_2)u_s$ and $v_s\in \big[-\sigma_2\eta,\sigma_1\eta\big]$, $s\in[0,\infty)$ with $\eta=\bar\mu +\delta(\sigma_1+\sigma_2)>0$. 
Observe that also in the extended model $Z^{(1)}$ does not depend on the control $v=(v_s)_{s\in[0,\infty)}$ and $Z^{(2),v}$ does not depend on the Brownian motion $(W_t)_{t\in[0,\infty)}$.

Using the same arguments as in Section \ref{sec:rho1_simple} (but with more lengthy computations) we obtain the following result.

\begin{theo}\label{value_function_extended}\hfill
	Let $y>\frac{\sigma_2}{\sigma_1}\,x$. Then the value function  for the extended model \eqref{model_mod}  is given by
		\begin{align*}
			V(x,y;\sigma_1,\sigma_2)=  \Phi&\left(\frac{\delta\sigma_2\,x
	+(\bar\mu+\delta\sigma_2)\,y}{N}\right)
-\exp\left(-\frac{2\bar\mu\,(x+y)}{(\sigma_1+\sigma_2)^2}\right)\Phi\left(\frac{A\sigma_2\,x+B\,y}{(\sigma_1+\sigma_2)\,N}\right)\\[0.3cm]
&-\exp\left(-\frac{2(\bar\mu+\delta\sigma_2)\,x}{\sigma_1^2}\right)\Phi\left(\frac{(\bar\mu+\delta\sigma_2)\,y-(A+\delta\sigma_2)\,\frac{\sigma_2}{\sigma_1}\,x}{N}\right)\\[0.3cm]
&+\exp\left(-\frac{2\bar\mu\, y}{(\sigma_1+\sigma_2)^2} -\frac{2\sigma_2\,x}{\sigma_1^2}\left(\frac{\bar\mu\,\sigma_2}{(\sigma_1+\sigma_2)^2}+\delta\right)\right)\Phi\left(\frac{B\,y-C\,\frac{\sigma_2}{\sigma_1}\,x}{(\sigma_1+\sigma_2)\,N}\right), 
	\end{align*}
		where 
\begin{align*}
	N&=\sqrt{\sigma_2(\sigma_1 y-\sigma_2 x)\big(\bar\mu+\delta(\sigma_1+\sigma_2)\big)}, \\[0.2cm]
	A&=2\bar\mu+\delta(\sigma_1+\sigma_2),\\[0.2cm]
	B&=(\bar\mu+\delta\sigma_2)(\sigma_1+\sigma_2)-2\bar\mu\sigma_1,\\[0.2cm]
		C&=\delta \sigma_1^2+3\delta\sigma_1\sigma_2+2(\bar\mu+\delta\sigma_2)\sigma_2.
	\end{align*}
	 For $y=\frac{\sigma_2}{\sigma_1}\,x$ it holds that 
	\begin{align*}
	V\left(x,\frac{\sigma_2}{\sigma_1}\,x; \sigma_1, \sigma_2\right)=1-\exp\left(-\frac{2\,\bar\mu\, x}{\sigma_1(\sigma_1+\sigma_2)}\right)= 1-\exp\left(-\frac{2\,\bar\mu\, y}{\sigma_2(\sigma_1+\sigma_2)}\right).
	\end{align*}
	For $y<\frac{\sigma_2}{\sigma_1}\,x$ we have
	\begin{align*}
	V(x,y;\sigma_1, \sigma_2)=V(y,x; \sigma_2, \sigma_1).
	\end{align*}
The optimal strategy for $V(x,y;\sigma_1,\sigma_2)$ is given by 
\begin{align*}
u^*_s=u^*_s\left(\widetilde{X}^{x,u^*}_s,\widetilde{Y}^{y,u^*}_s\right)= \begin{cases}
\bar\mu+\delta\sigma_2, \qquad&\widetilde{Y}^{y,u^*}_s\geq \frac{\sigma_2}{\sigma_1}\,\widetilde{X}^{x,u^*}_s,\\[0.3cm]
-\delta\sigma_1, &\widetilde{Y}^{y,u^*}_s< \frac{\sigma_2}{\sigma_1}\,\widetilde{X}^{x,u^*}_s.
\end{cases}
\end{align*}

\end{theo}

{\captionsetup[figure]{skip=4pt}
	\begin{figure}
		\begin{subfigure}{0.485\textwidth}
			\includegraphics[width=\textwidth]{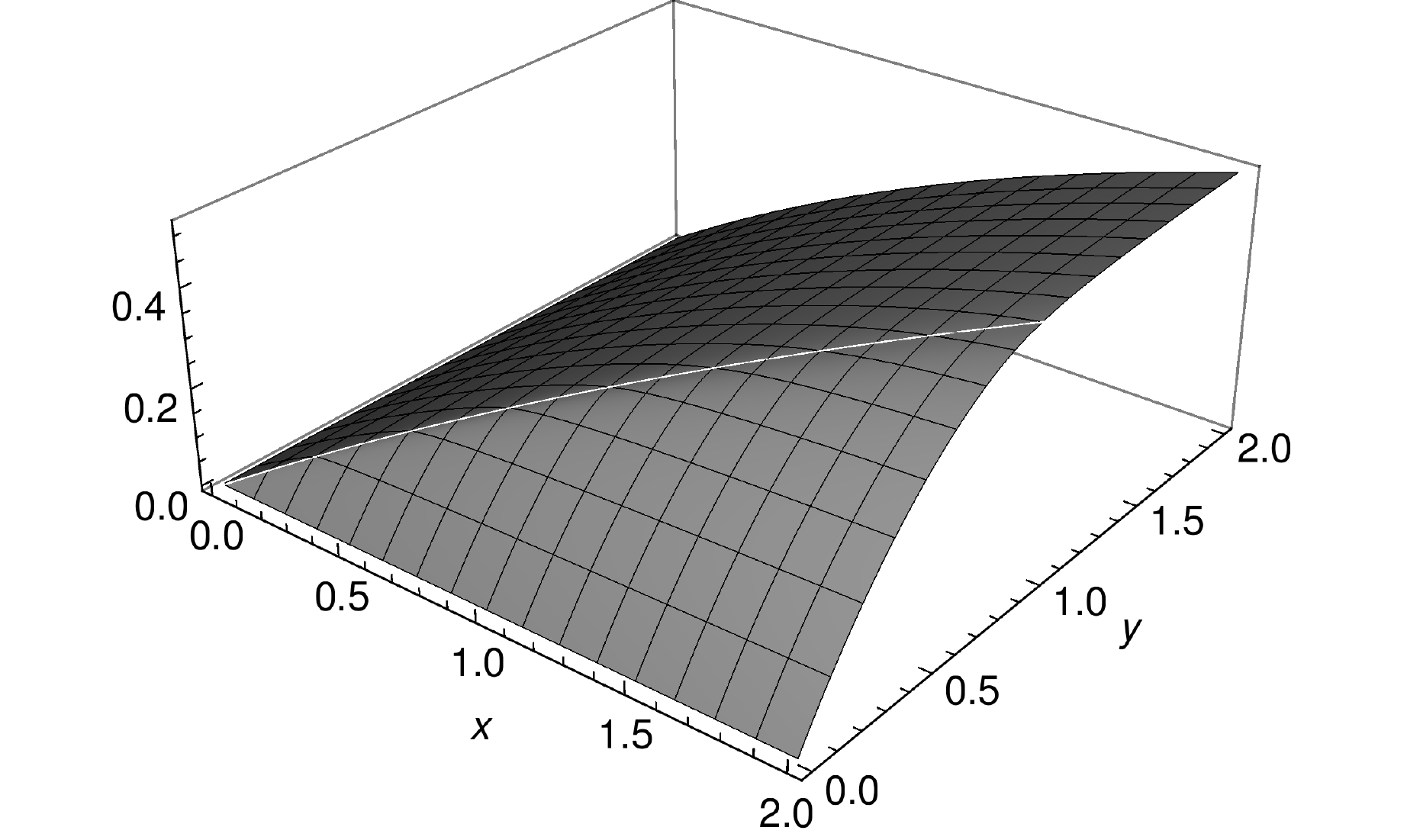}
			\subcaption{$\sigma_1=2,\,\sigma_2=1,\, \delta=-\frac 14$.}
		\end{subfigure}
		\hspace{0.02\textwidth}
		\begin{subfigure}{0.485\textwidth}
			\includegraphics[width=\textwidth]{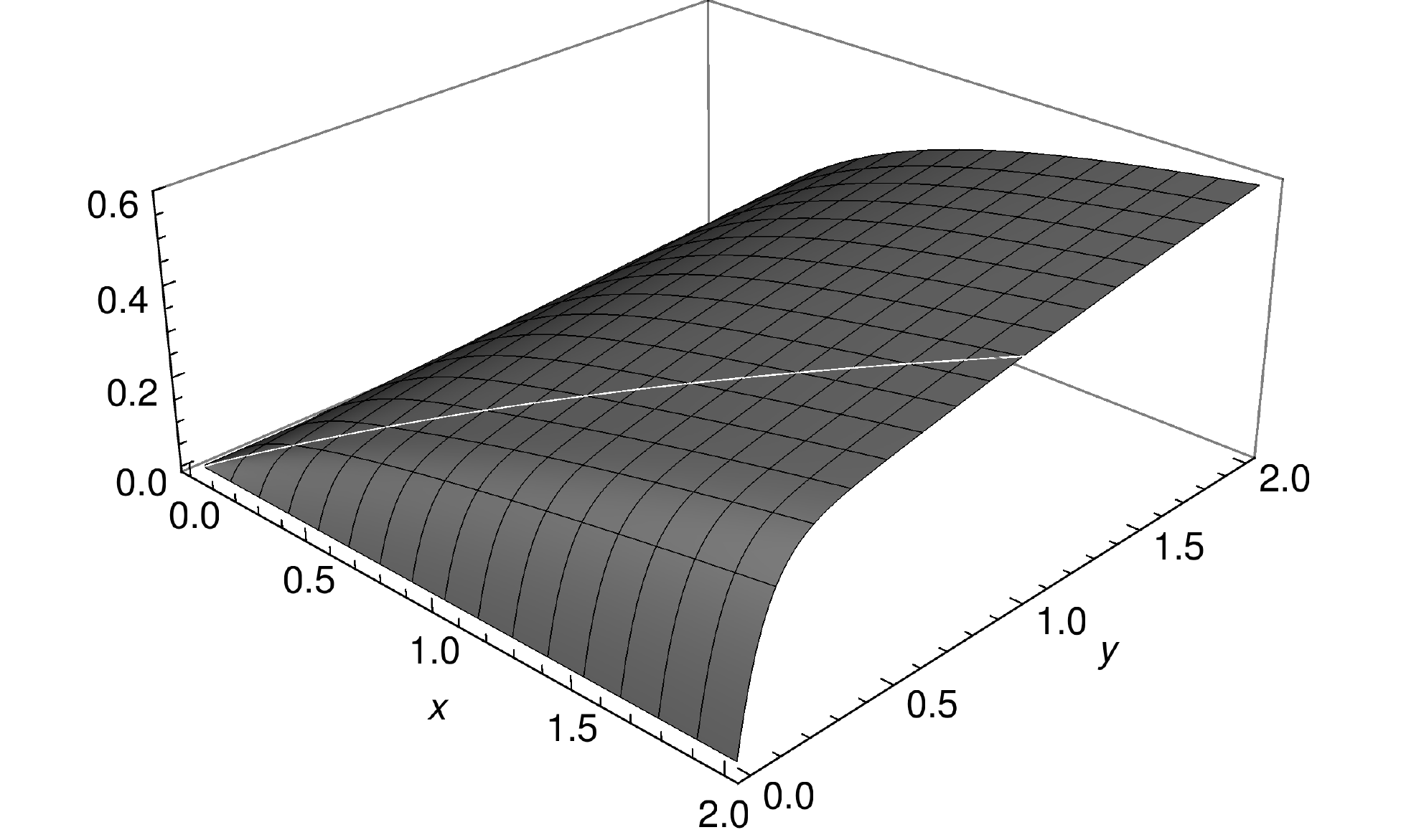}
			\subcaption{$\sigma_1=2,\, \sigma_2=1,\,\delta=2$.}
		\end{subfigure}
		
		\begin{subfigure}{0.485\textwidth}
			\includegraphics[width=\textwidth]{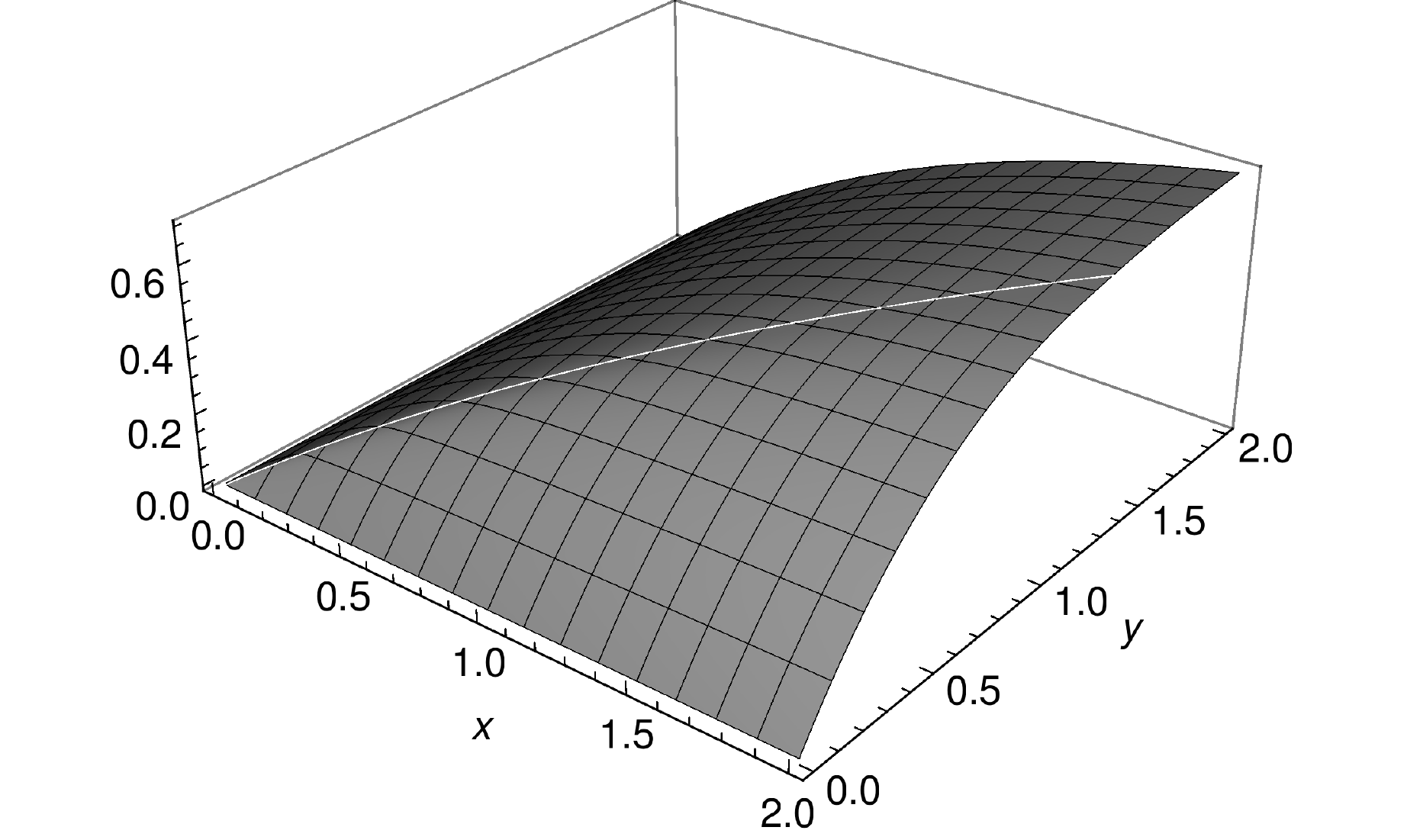}
			\subcaption{$\sigma_1=\frac32,\, \sigma_2=1,\, \delta=-\frac 14$.}
		\end{subfigure}
		\hspace{0.02\textwidth}
		\begin{subfigure}{0.485\textwidth}
			\includegraphics[width=\textwidth]{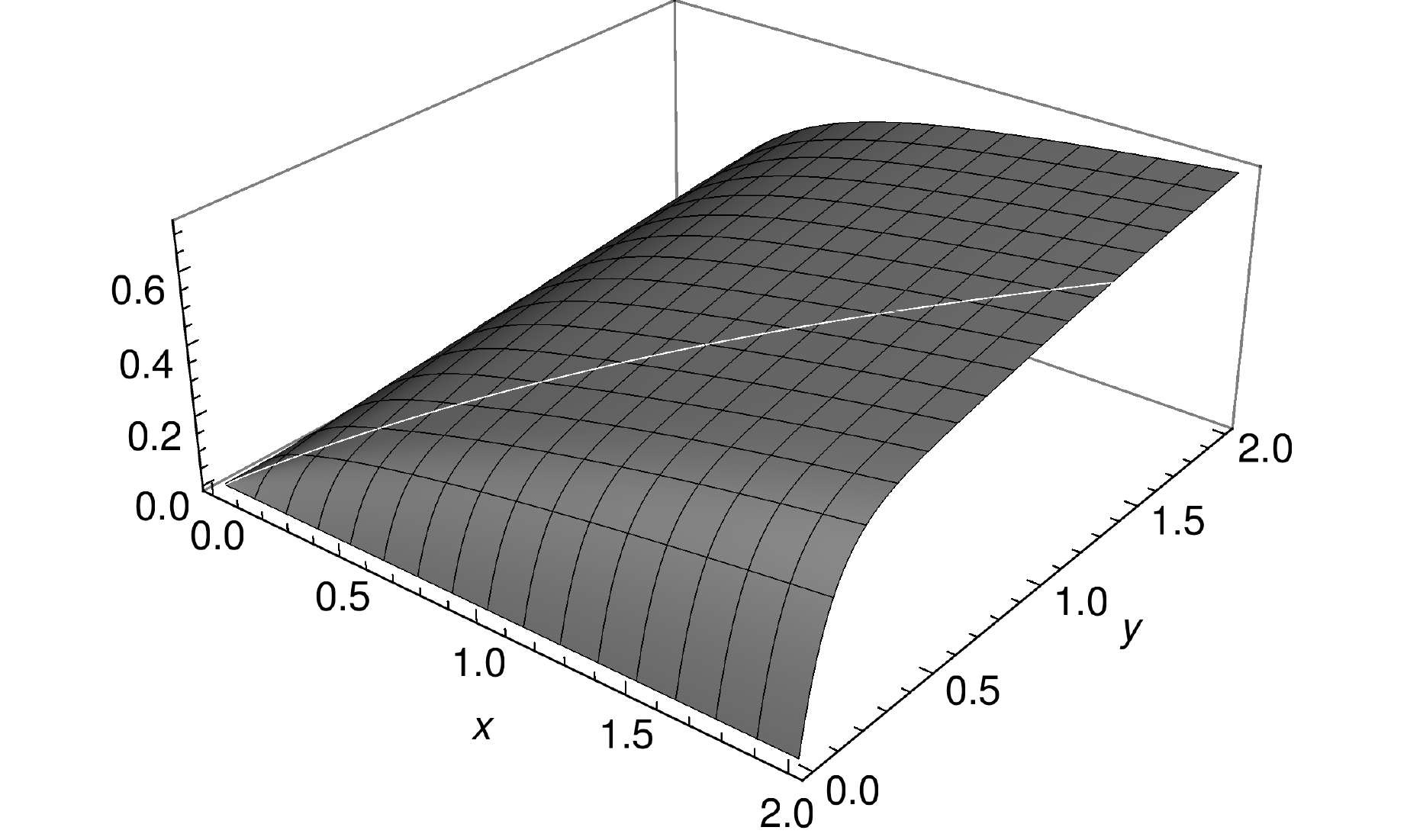}
			\subcaption{$\sigma_1=\frac 32,\,\sigma_2=1,\, \delta=2$.}
		\end{subfigure}
		\caption{The value function for $\rho=1$, $\bar\mu=1$, different diffusion rates $\sigma_1, \sigma_2>0$ and different $\delta>-\frac{\bar\mu}{\sigma_1+\sigma_2}$. }\label{fig:value_mod}
	\end{figure}
\begin{remark}
	Observe that the Hamilton-Jacobi-Bellman equation in the extended model \eqref{model_mod} is given by 
\begin{align}\label{HJB_mod}
\left(\frac{\sigma_1^2}{2} V_{xx}+\frac{\sigma_2^2}{2} V_{yy}+\sigma_1\sigma_2 V_{xy}+\max_{a\in [- \delta\sigma_1,\, \bar\mu+\delta\sigma_2]} \left\{a V_x+ (\bar\mu-a) V_y\right\}\right)(x,y;\sigma_1,\sigma_2)&=0
\end{align}
for $x,y>0$ with boundary conditions
\begin{align}
V(x,0; \sigma_1,\sigma_2)=V(0,y;\sigma_1,\sigma_2)&=0, \label{boundary1}\\
\lim_{x\to\infty} V(x,y;\sigma_1,\sigma_2)&= 1-\exp\left(-\frac{2(\bar\mu+\delta\sigma_1)\,y}{\sigma_2^2}\right),\label{boundary2} \\
\lim_{y\to\infty} V(x,y;\sigma_1,\sigma_2)&= 1-\exp\left(-\frac{2(\bar\mu+\delta\sigma_2)\,x}{\sigma_1^2}\right).\label{boundary3}
\end{align}
One can prove that $V(\cdot\,,\cdot\,;\sigma_1,\sigma_2)\in C^2\big((0,\infty)^2\big)\,\cap\, C\big([0,\infty)^2\big)$ and that $V(\cdot\,,\cdot\,;\sigma_1,\sigma_2)$ solves~\eqref{HJB_mod} with boundary conditions \eqref{boundary1}, \eqref{boundary2} and \eqref{boundary3}. Furthermore,  
\begin{align}
\left\{(x,y)\in(0,\infty)^2\colon V_x(x,y;\sigma_1, \sigma_2)> V_y(x,y;\sigma_1, \sigma_2)\right\}
&=\left\{(x,y)\in(0,\infty)^2\colon y>\frac{\sigma_2}{\sigma_1}\,x\right\},\label{derivatives} \\
\left\{(x,y)\in(0,\infty)^2\colon V_x(x,y;\sigma_1, \sigma_2)< V_y(x,y;\sigma_1, \sigma_2)\right\}
&=\left\{(x,y)\in(0,\infty)^2\colon y<\frac{\sigma_2}{\sigma_1}\,x\right\},\label{derivatives2} \\
\left\{(x,y)\in(0,\infty)^2\colon V_x(x,y;\sigma_1, \sigma_2)=V_y(x,y;\sigma_1, \sigma_2)\right\}
&=\left\{(x,y)\in(0,\infty)^2\colon y=\frac{\sigma_2}{\sigma_1}\,x\right\}.\notag  
\end{align}
To see that \eqref{derivatives} holds true, note that for $y>\frac{\sigma_2}{\sigma_1}x$ we have
\begin{align*}
V_x&(x,y; \sigma_1, \sigma_2)-V_y(x,y; \sigma_1, \sigma_2)\\[0.3cm]
=&\,\frac{2(\bar\mu+\delta\sigma_2)}{\sigma_1^2}\exp\left(-\frac{2(\bar\mu+\delta\sigma_2)\,x}{\sigma_1^2}\right)\Phi\left(\frac{(\bar\mu+\delta\sigma_2)\,y-(A+\delta\sigma_2)\frac{\sigma_2}{\sigma_1}\,x}{N}\right)\\[0.3cm]
& +\frac{2\big(\bar\mu(\sigma_1-\sigma_2)-\delta\sigma_2(\sigma_1+\sigma_2)\big)}{\sigma_1^2(\sigma_1+\sigma_2)}\exp\left(-\frac{2\bar\mu\, y}{(\sigma_1+\sigma_2)^2} -\frac{2\sigma_2\,x}{\sigma_1^2}\left(\frac{\bar\mu\,\sigma_2}{(\sigma_1+\sigma_2)^2}+\delta\right)\right)\Phi\left(\frac{ B\, y-C\,\frac{\sigma_2}{\sigma_1}\,x}{(\sigma_1+\sigma_2)\,N}\right). 
\end{align*}
If $\delta\in \left(-\frac{\bar\mu}{\sigma_1+\sigma_2},\frac{\bar\mu(\sigma_1-\sigma_2)}{\sigma_2(\sigma_1+\sigma_2)} \right]$, then one can directly conclude that $(V_x-V_y)(x,y;\sigma_1,\sigma_2)>0$ for all $y>\frac{\sigma_2}{\sigma_1}\,x$. 
If $\delta>\frac{\bar\mu(\sigma_1-\sigma_2)}{\sigma_2(\sigma_1+\sigma_2)} $, then observe that 
\begin{align*}
\lim_{y\downarrow \frac{\sigma_2}{\sigma_1}x} (V_x-V_y)(x,y;\sigma_1,\sigma_2)=0
\end{align*}
and
\begin{align*}
(&V_x-V_y)_y(x,y;\sigma_1, \sigma_2)\\[0.3cm]
&= \frac{4\bar\mu\left(\delta\sigma_2(\sigma_1+\sigma_2)+\bar\mu(\sigma_2-\sigma_1)\right)}{\sigma_1^2(\sigma_1+\sigma_2)^3}\exp\left(-\frac{2\bar\mu\, y}{(\sigma_1+\sigma_2)^2} -\frac{2\sigma_2\,x}{\sigma_1^2}\left(\frac{\bar\mu\,\sigma_2}{(\sigma_1+\sigma_2)^2}+\delta\right)\right)\Phi\left(\frac{ B\, y-C\,\frac{\sigma_2}{\sigma_1}\,x}{(\sigma_1+\sigma_1)\,N}\right)\\[0.3cm]
&\quad \ +\sqrt{\frac{2}{\pi}}
\big(2\sigma_1\, y+(\sigma_1-\sigma_2)\,x\big)\frac{\bar\mu\sqrt{\sigma_2(\bar\mu+\delta(\sigma_1+\sigma_2))}}{\sigma_1(\sigma_1+\sigma_2)^2(\sigma_1\, y-\sigma_2\,x)^{\frac 32}}\exp\left(-\frac{\big(\delta \sigma_2\, x+(\bar\mu+\delta\sigma_2)\, y\big)^2}{2N^2}\right) >0 
\end{align*}
for $y>\frac{\sigma_2}{\sigma_1}\,x$. Hence, also in this case we have $(V_x-V_y)(x,y;\sigma_1,\sigma_2)>0$.

Similarly, one derives \eqref{derivatives2}. 
\end{remark}

Figure \ref{fig:value_mod} depicts the value function $V(x,y;\sigma_1, \sigma_2)$ for different diffusion rates and different $\delta$.

 If we consider two insurance companies whose endowment processes are Brownian motions with drift  $\mu_i>0$, diffusion coefficients $\sigma_i>0$, $i=1,2$, and initial endowments $x,y>0$, respectively,  
then the survival probability can be derived similarly to the case $\sigma_1=\sigma_2=1$ and is given by
\begin{align*}
\P[&x+\mu_1\, t +\sigma_1 W_t>0,\, y+\mu_2\, t+\sigma_2 W_t \text{ for all $t\in[ 0,\infty)$}]\\[0.2cm]
&=\begin{cases}
\begin{aligned}
&1-\exp\left(-2\min\left\{\frac{\mu_1\,x }{\sigma_1^2}, \frac{\mu_2\, y}{\sigma_2^2}\right\}\right), && \frac{\mu_1}{\sigma_1}=\frac{\mu_2}{\sigma_2}; \\[0.2cm]
& && \frac{\mu_1}{\sigma_1}>\frac{\mu_2}{\sigma_2},   \frac{x}{\sigma_1}\geq\frac{y}{\sigma_2};\\[0.2cm]
& && \frac{\mu_1}{\sigma_1}<\frac{\mu_2}{\sigma_2}, \frac{x}{\sigma_1}\leq\frac{y}{\sigma_2};\\[0.5cm]
&\Phi\left(\frac{|\mu_1\, y-\mu_2\, x|}{L}\right)-\exp\left(-\frac{2\mu_1 \,x}{\sigma_1^2}\right)\Phi\left(\frac{\mu_1\, y +D(2,1)\,x}{L}\right)\qquad\qquad  &&\frac{\mu_1}{\sigma_1}>\frac{\mu_2}{\sigma_2}, \frac{x}{\sigma_1}<\frac{y}{\sigma_2};\\[0.2cm]
&\quad -\exp\left(-\frac{2\mu_2 \,y}{\sigma_2^2}\right)\Phi\left(\frac{\mu_2 \,x +D(1,2)\,y}{L}\right) &&\frac{\mu_1}{\sigma_1}<\frac{\mu_2}{\sigma_2}, \frac{x}{\sigma_1}>\frac{y}{\sigma_2};\\[0.3cm]
&\quad+\exp\left(-\frac{2\mu_1\, x}{\sigma_1^2}-\frac{2\mu_2\, y}{\sigma_2^2}+\frac{ 4\min\left\{\mu_2\, x, \mu_1\,y\right\}}{\sigma_1\sigma_2}\right)&&\hspace*{-3.1cm}\Phi\left(\frac{D(2,1)\,x+D(1,2)\,y+2\min\left\{\mu_2\, x , \mu_1\,y \right\}}{L}\right),
&
\end{aligned}
\end{cases}
\end{align*}
where 
\begin{align*}
L=\sqrt{(\sigma_1\,y-\sigma_2\, x)(\mu_1\sigma_2-\mu_2\sigma_1)},\qquad 
D(i,j)=\mu_i-\frac{2\mu_j\sigma_i}{\sigma_j},\quad  i,j\in\{1,2\}.
\end{align*}
The gain of collaboration is then given by the difference 
\begin{align*}
V(x,y;\sigma_1,\sigma_2)-\P[x+\mu_1\, t +\sigma_1 W_t>0, y+\mu_2\, t+\sigma_2 W_t \,\text{ for all $t\in[0,\infty)$}].
\end{align*}
See Figure \ref{fig:gain_mod} for the gain of collaboration in case of different diffusion coefficients and drift rates. 

\begin{remark}
		Observe that in the extended model we also relax the assumption on $\delta$ and do not assume that $\delta>-\min\left\{\frac{\mu_1}{\sigma_1},\frac{\mu_2}{\sigma_2}\right\}$. For $\delta\in \left(-\frac{\bar\mu}{\sigma_1+\sigma_2}, -\min\left\{\frac{\mu_1}{\sigma_1}, \frac{\mu_2}{\sigma_2}\right\}\right)$ the companies are forced to collaborate, in particular firm $i$ with $\frac{\mu_i}{\sigma_i}>\frac{\mu_j}{\sigma_j}$, $j\neq i$, has to make payments to firm $j$ at every point in time.  Thus, the gain of collaboration can be negative, see Figure \ref{fig:neg_gain}. 
\end{remark}

{\captionsetup[figure]{skip=4pt}
	\begin{figure}
		\begin{subfigure}{0.485\textwidth}
		\includegraphics[width=\textwidth]{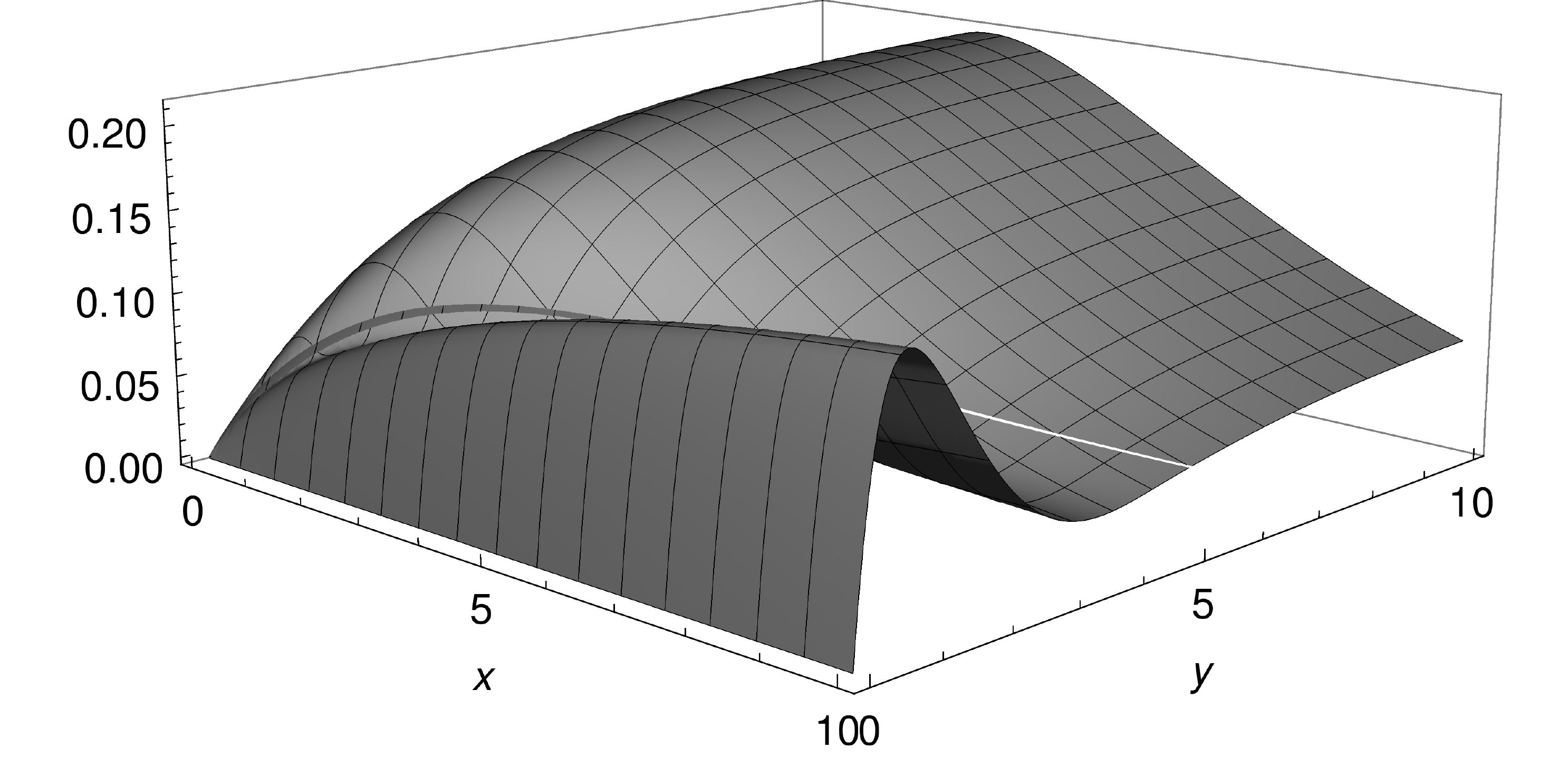}
			\subcaption{$\mu_1=\mu_2=\frac 12,\,\sigma_1=2,\, \sigma_2=1,\, \delta=-\frac {1}{10} $.}
		\end{subfigure}
		\hspace{0.02\textwidth}
		\begin{subfigure}{0.485\textwidth}
		\includegraphics[width=\textwidth]{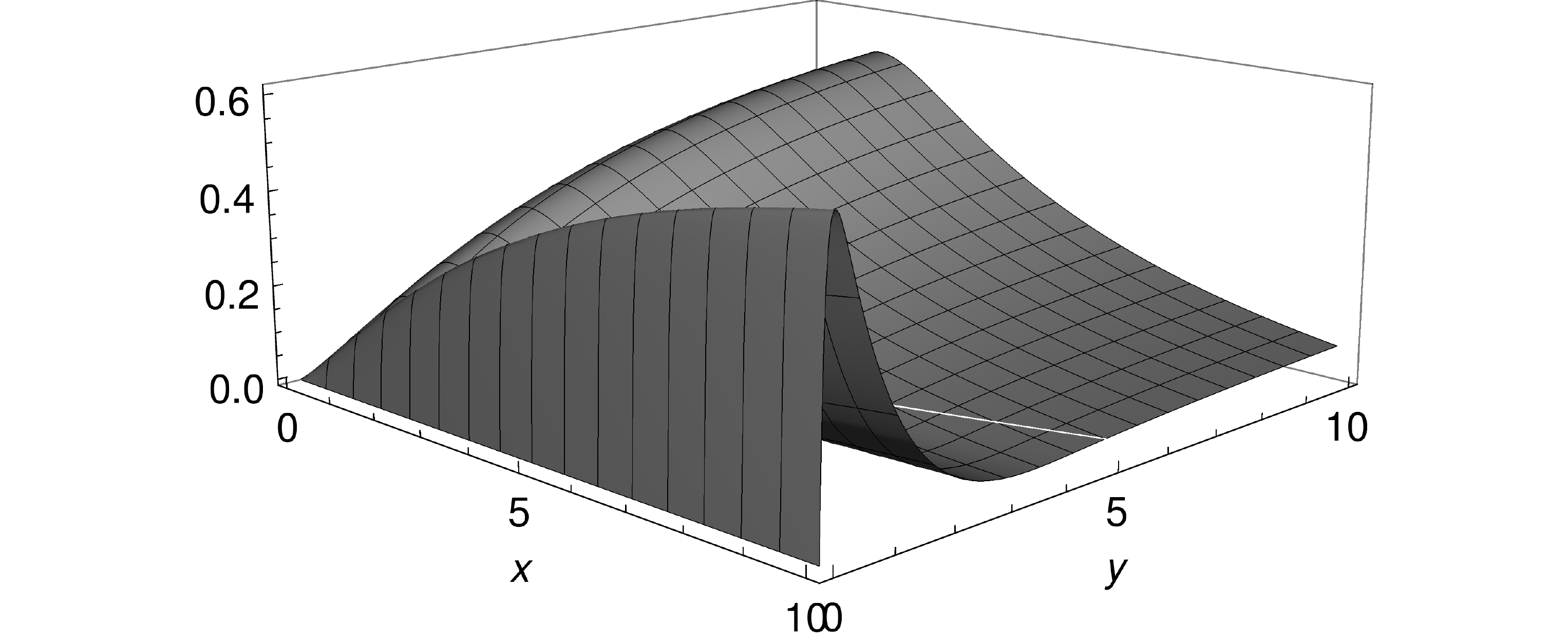}
		\vspace{0.05cm}
			\subcaption{$\mu_1=\mu_2=\frac 12,\,\sigma_1=2,\, \sigma_2=1,\, \delta=2$.}
		\end{subfigure}
		\\[0.2cm]
		\begin{subfigure}{0.49\textwidth}
			\includegraphics[width=\textwidth]{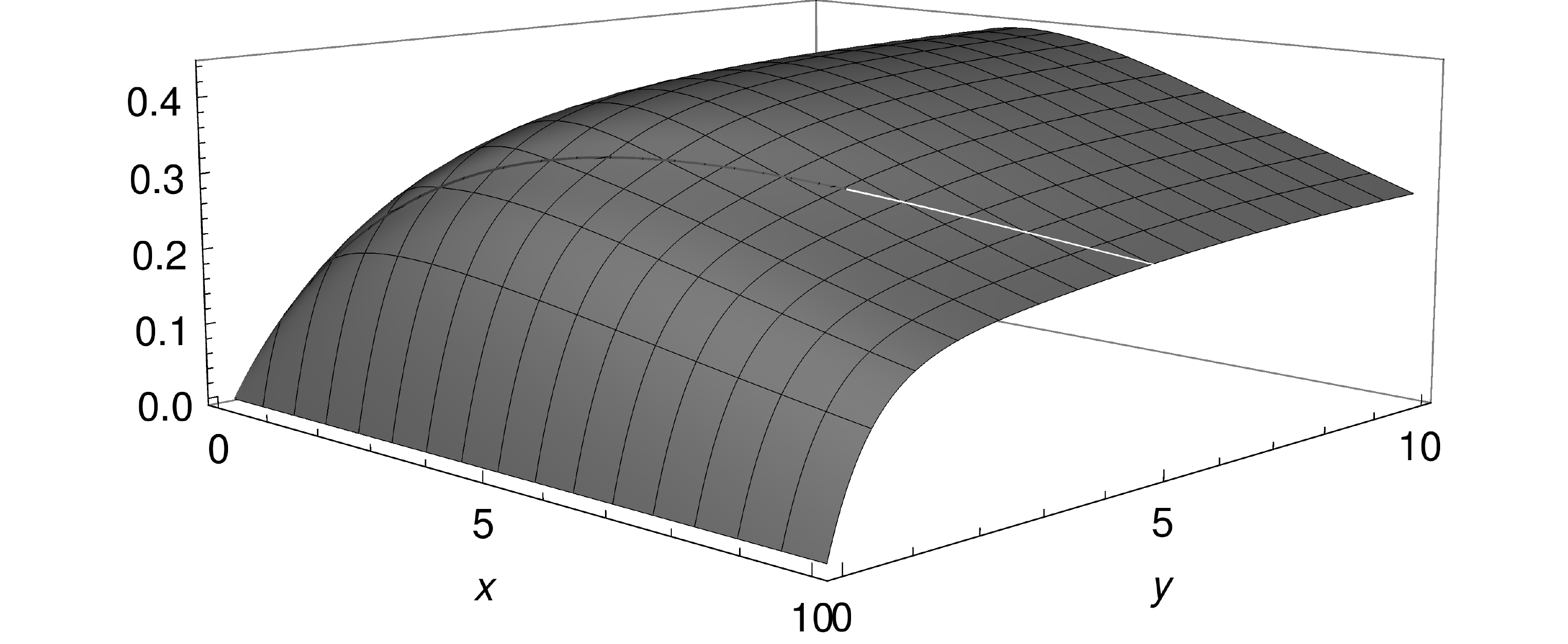}
				\subcaption{$\mu_1=\frac 14,\,\mu_2=\frac 34,\sigma_1=2,\, \sigma_2=1,\,\delta=-\frac{1}{10}$.}
		\end{subfigure}
		\hspace{0.02\textwidth}
		\begin{subfigure}{0.485\textwidth}
				\includegraphics[width=\textwidth]{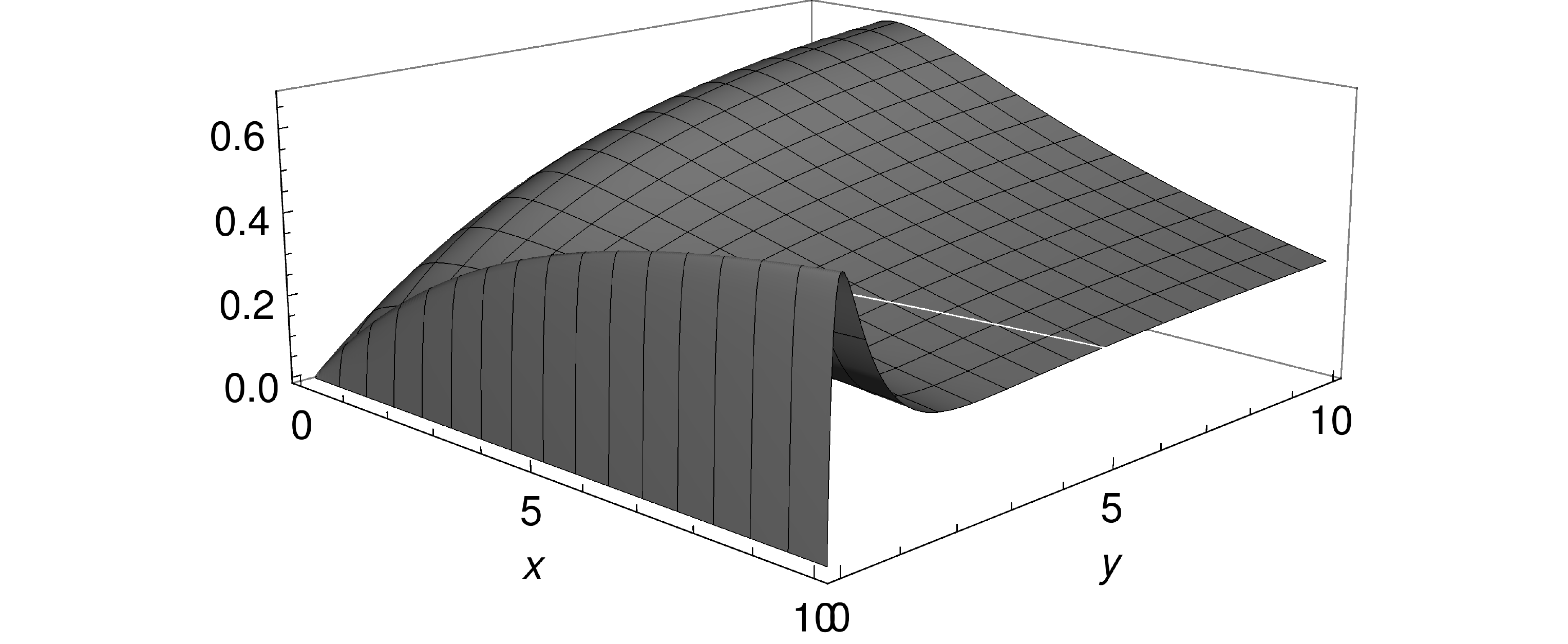}
		\vspace{-0.45cm}
		\subcaption{$\mu_1=\frac 14,\, \mu_2=\frac 34,\,\sigma_1=2,\,\sigma_2=1,\, \delta=2$.}
		\end{subfigure}
		\\[0.2cm]
		\begin{subfigure}{0.485\textwidth}
			
			\includegraphics[width=\textwidth]{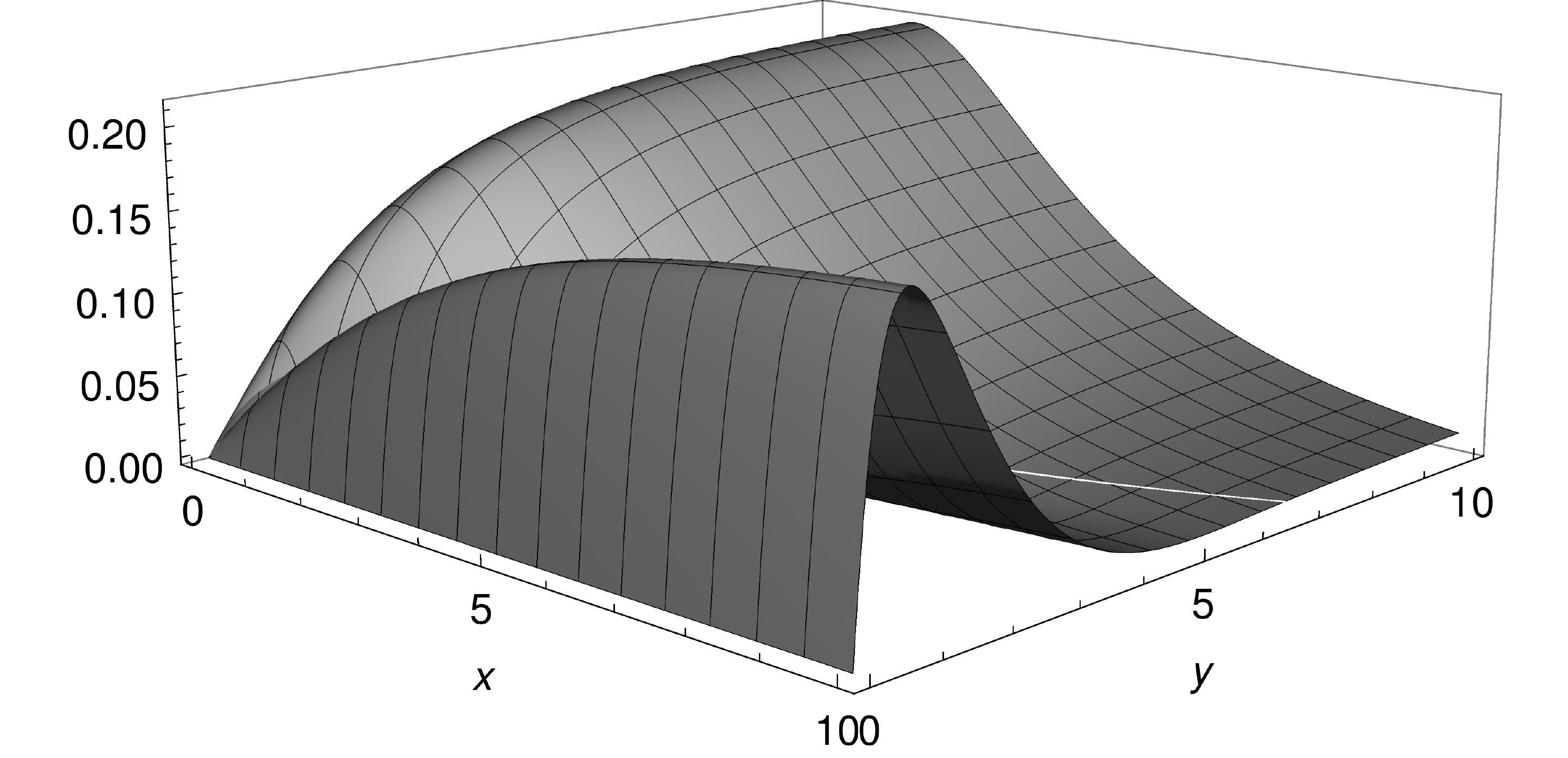}
	
			\subcaption{$\mu_1=\mu_2=\frac 12,\,\sigma_1=\frac 32,\, \sigma_2=1,\, \delta=-\frac{1}{10} $.}
		\end{subfigure}
		\hspace{0.02\textwidth}
		\begin{subfigure}{0.485\textwidth}
			\includegraphics[width=\textwidth]{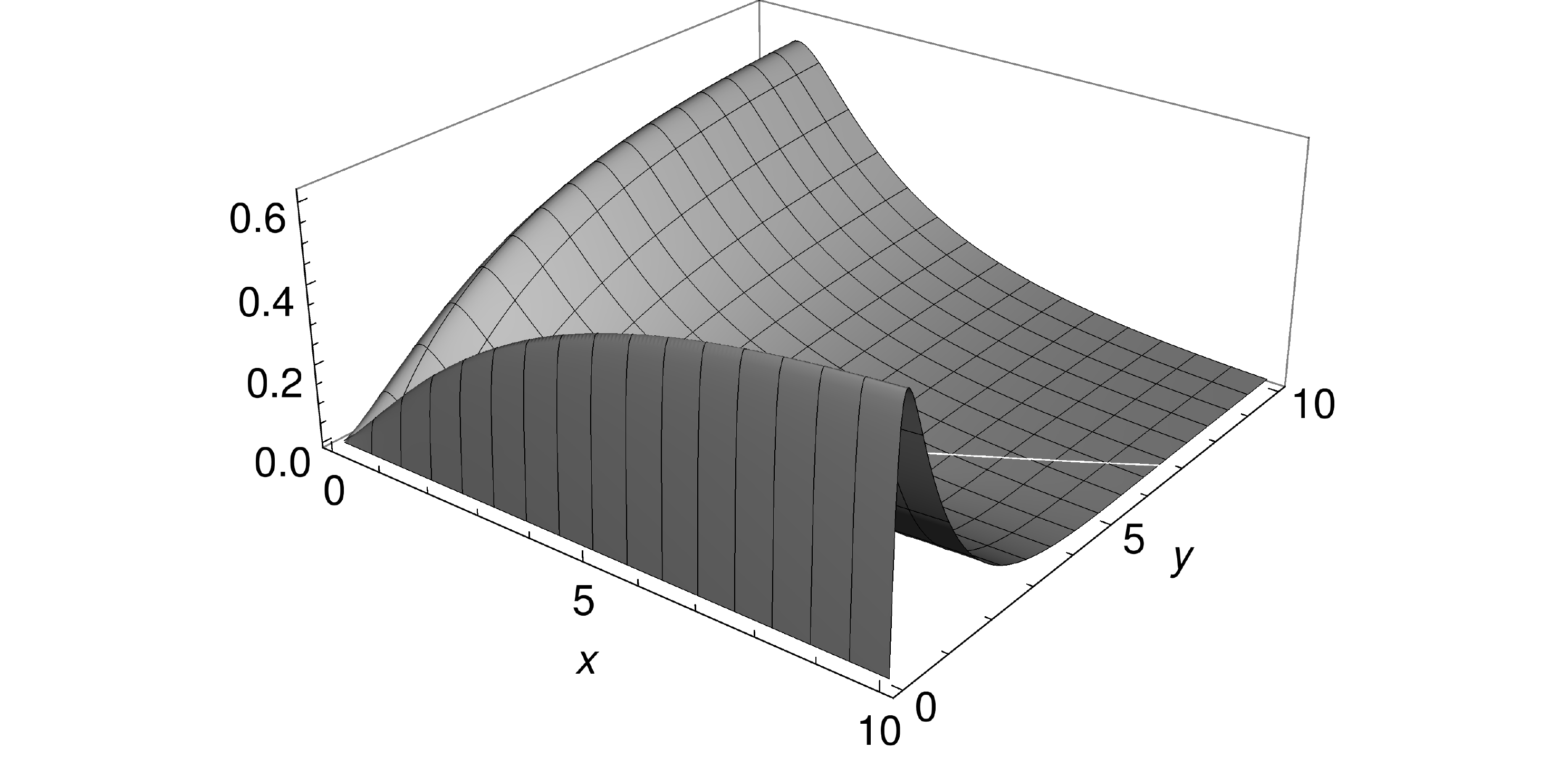}
			\subcaption{$\mu_1=\mu_2=\frac 12,\,\sigma_1=\frac32,\, \sigma_2=1,\, \delta=2$.}
		\end{subfigure}
		\\[0.2cm]
		\begin{subfigure}{0.485\textwidth}
		
			\includegraphics[width=\textwidth]{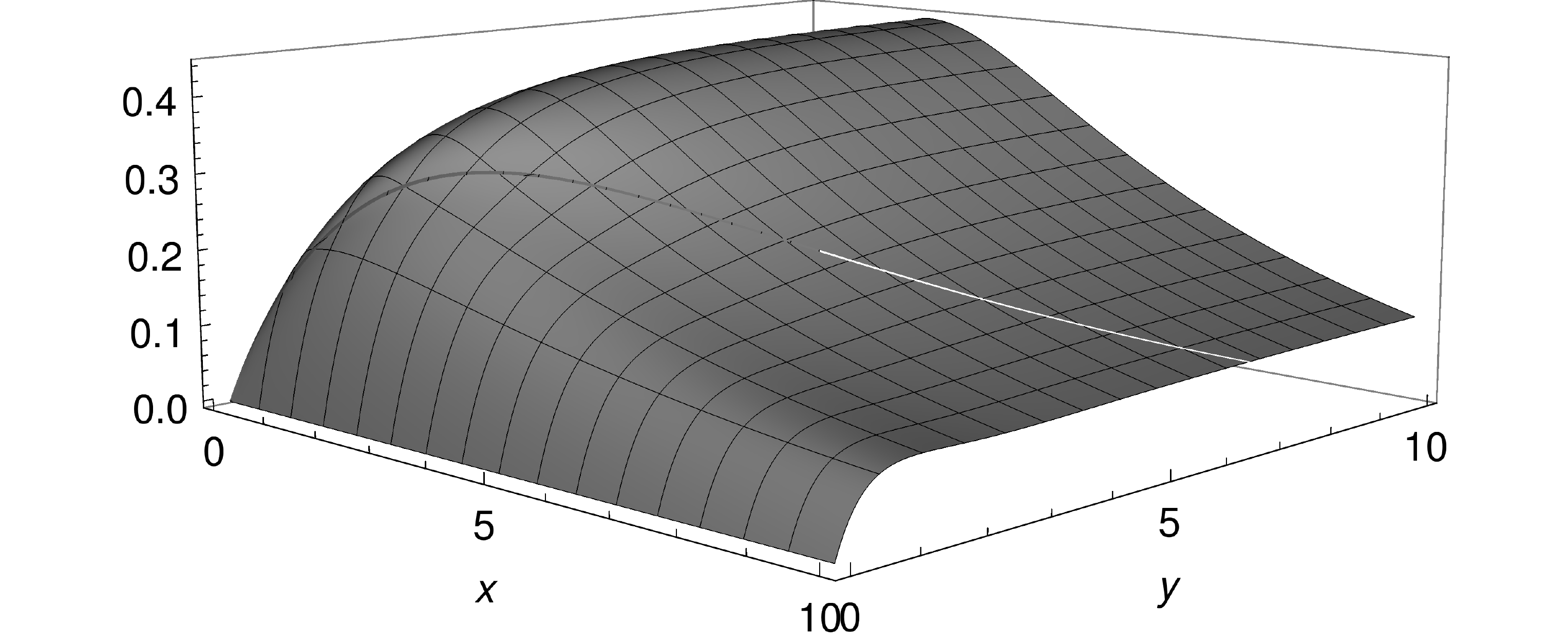}
	
			\subcaption{$\mu_1=\frac 14,\,\mu_2=\frac 34,\,\sigma_1=\frac 32,\, \sigma_2=1,\,\delta=-\frac{1}{10}$.}
		\end{subfigure}
		\hspace{0.02\textwidth}
		\begin{subfigure}{0.485\textwidth}
			\includegraphics[width=\textwidth]{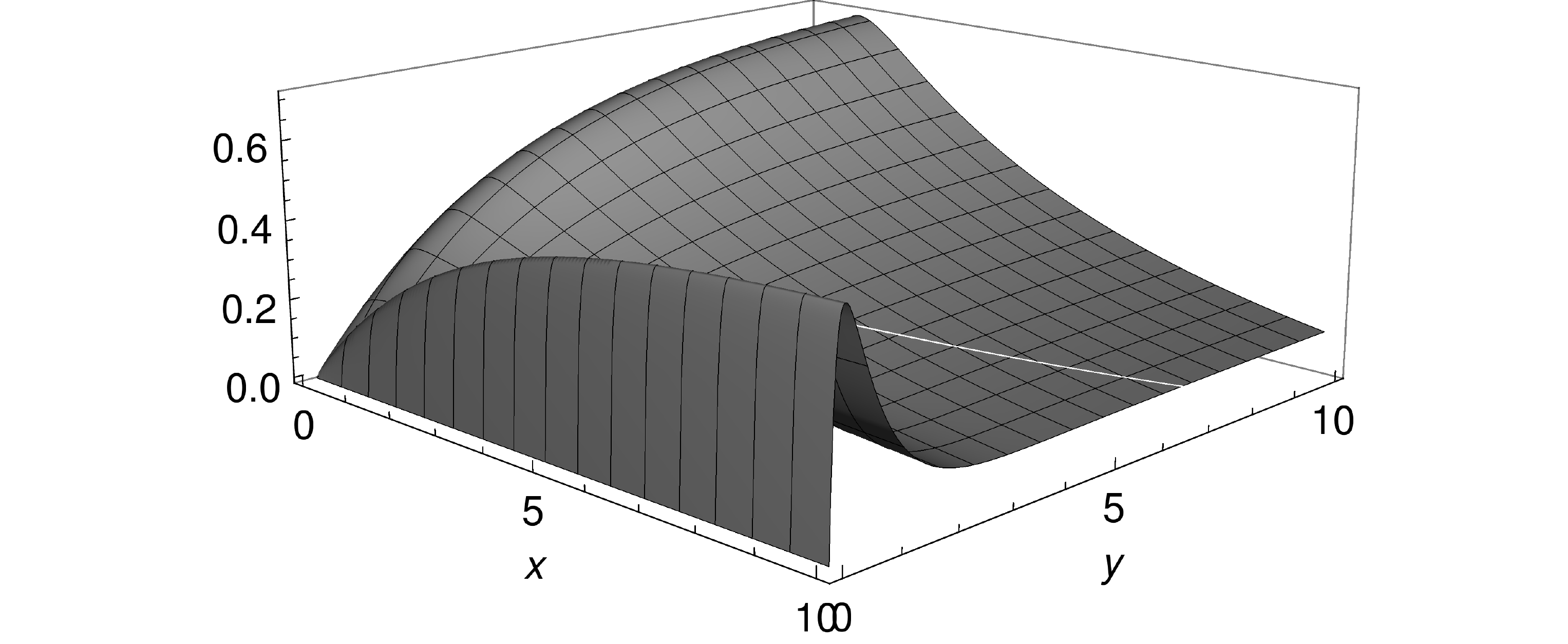}
			\subcaption{$\mu_1=\frac 14,\, \mu_2=\frac 34,\,\sigma_1=\frac 32,\,\sigma_2=1,\, \delta=2$.}
		\end{subfigure}
		\caption{The gain of collaboration for $\rho=1$, $\bar\mu=1$ and different drift rates $\mu_1$, $\mu_2>0$, diffusion rates $\sigma_1, \sigma_2>0$ and different $\delta>-\frac{\bar\mu}{\sigma_1+\sigma_2}$. \vspace*{0.07cm}}\label{fig:gain_mod}
	\end{figure}
}

\begin{figure}
	\begin{center}
	\includegraphics[width=0.6\textwidth]{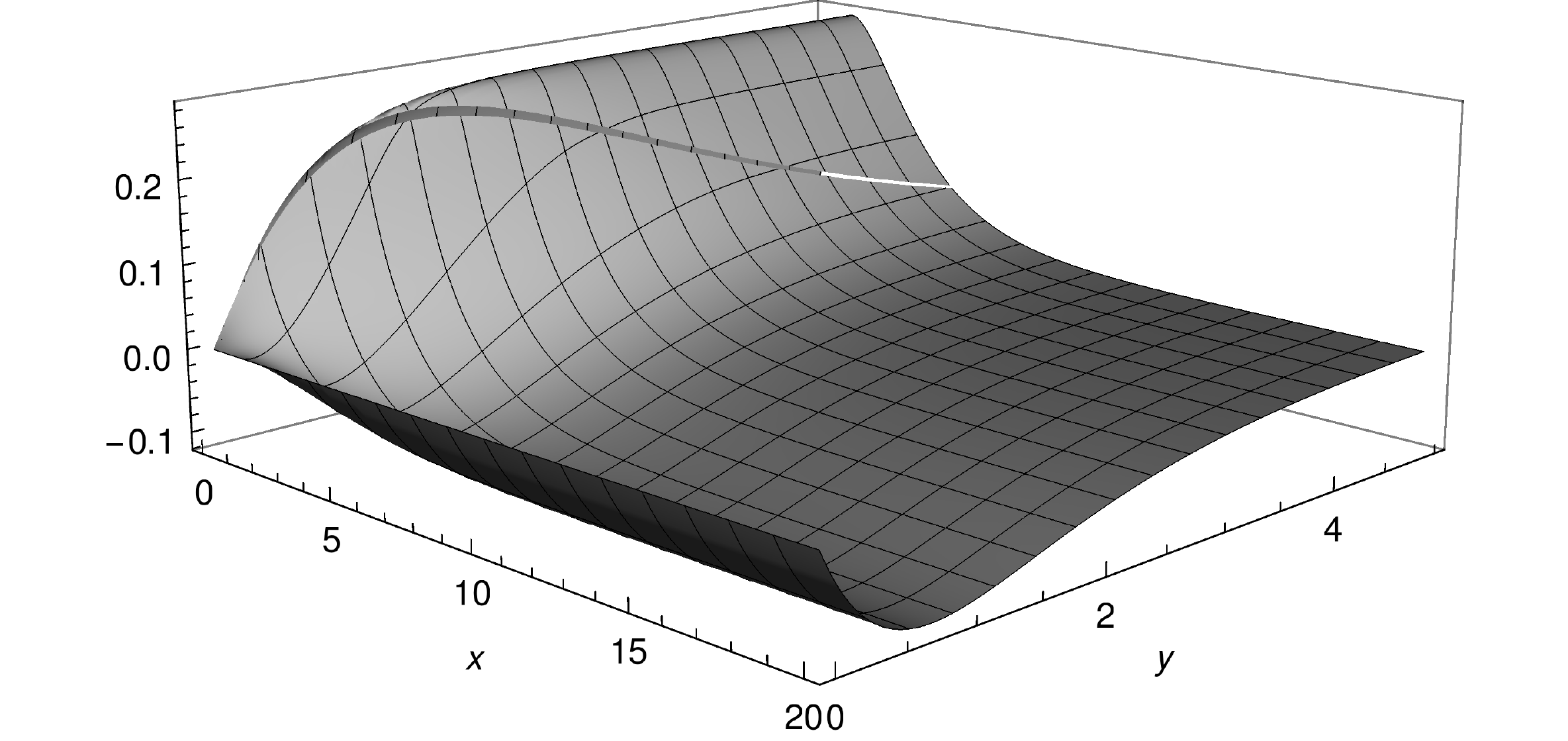}	
	\caption{For $\mu_1=\frac 14$, $\mu_2=\frac34$, $\sigma_1=\sigma_2=1$, $\delta= -\frac{9}{20}$ the gain of collaboration can be negative.}\label{fig:neg_gain}
	\end{center}
\end{figure}

\section{Perfectly Negative Correlation: $\rho=-1$}\label{sec:rho-1}
In this section we focus on the case $\rho=-1$ and obtain a different characterization of the value function in terms of the probability that a reflected Brownian motion with drift never hits a specific line. Unfortunately, we cannot use similar arguments as in the case $\rho=1$ to derive an explicit formula for the value function. 

From Theorem \ref{optimal} we already know that an optimal strategy for the transfer payments is given by
\begin{align*}
u^*_s=(\bar\mu+\delta)\,\mathds{1}_{\left\{X_s^{x,u^*}\!\leq\,Y^{y,u^*}_s\right\}}-\delta\, \mathds{1}_{\left\{X_s^{x,u^*}\!>\,Y^{y,u^*}_s\right\}}.
\end{align*}
 For perfectly negatively correlated Brownian motions and for the optimal strategy $u^*$ it holds that 
\begin{align*}
Z^{(1)}_t&= X_t^{x,u^*}+Y_t^{y,u^*}= x+y+\bar\mu\, t , \\[0.3cm]
Z^{(2)}_t&= Y_t^{y,u^*}-X_t^{x,u^*}= y-x-2W_t+\int_0^t (\bar\mu-2u_s^*)\,ds= y-x-2W_t-\int_0^t (\bar\mu+2\delta)\,\text{sign}\big(Z^{(2)}_t\big)\, ds.
\end{align*}
Moreover, 
\begin{align*}
\tau(x,y;u^*)
&= \inf\left\{t\in[0,\infty)\colon \left|Z^{(2)}_t\right|\geq x+y+\bar\mu\, t\right\}
\end{align*}
and $V(x,y)=\P[\tau(x,y;u^*)=\infty]$. 

The process $\big(\frac 12\big|Z^{(2)}_t\big|\big)_{t\in[0,\infty)}$ is a representation of a reflected Brownian motion with drift $-\left(\frac{\bar\mu}{2}+\delta\right)$, see~\cite{GraversenShiryaev}.  Therefore, the value function $V(x,y)$ can be interpreted as the probability that a reflected Brownian motion with drift  $-\left(\frac{\bar\mu}{2}+\delta\right)$, never hits the linear barrier 
\begin{align*}
b(t)=\frac 12(x+y+\bar\mu\,t).
\end{align*}

To the best of our knowledge no closed formula for the hitting probability is available in the literature.

\end{document}